\documentclass{amsart}
\usepackage{amssymb,mathrsfs,array,color}

\newcommand{\g}{\mathfrak{g}}
\newcommand{\h}{\mathfrak{h}}
\newcommand{\mfk}{\mathfrak{k}}
\newcommand{\mfp}{\mathfrak{p}}
\newcommand{\mft}{\mathfrak{t}}
\newcommand{\mfb}{\mathfrak{b}}
\newcommand{\mfn}{\mathfrak{n}}
\newcommand{\mfa}{\mathfrak{a}}

\newcommand{\bbZ}{\mathbb{Z}}
\newcommand{\bbC}{\mathbb{C}}
\newcommand{\bbA}{\mathbb{A}}
\newcommand{\bbQ}{\mathbb{Q}}
\newcommand{\bbR}{\mathbb{R}}
\newcommand{\cO}{\mathcal{O}}
\newcommand{\cF}{\mathcal{F}}
\newcommand{\ob}{\mathrm{Ob}}

\DeclareMathOperator{\Hom}{\mathrm{Hom}}
\DeclareMathOperator{\sgn}{\mathrm{sgn}}

\newtheorem{theorem}{Theorem}[section]
\newtheorem{proposition}[theorem]{Proposition}
\newtheorem{lemma}[theorem]{Lemma}
\newtheorem{corollary}[theorem]{Corollary}

\theoremstyle{definition}

\theoremstyle{remark}

\newtheorem{remark}[theorem]{Remark}

\title[Signatures of Invariant Hermitian Forms on Highest Weight Modules]
{Signatures of Invariant Hermitian Forms on Irreducible Highest Weight 
Modules}

\author{Wai Ling Yee}
\address{Department of Mathematical and Statistical Sciences \\
University of Alberta \\
Edmonton, Alberta, CANADA}
\email{wlyee@math.ualberta.ca}
\thanks{This research was supported by an NSERC postdoctoral fellowship and 
by the Killam Foundation.}

\begin{document}
\begin{abstract}
Perhaps the most important problem in representation theory in the 1970s
and early 1980s was the determination of the multiplicity of composition
factors in a Verma module.  This problem was settled by the proof of the
Kazhdan-Lusztig Conjecture which states that the multiplicities may be
computed via Kazhdan-Lusztig polynomials.  In this paper, we introduce
signed Kazhdan-Lusztig polynomials, a variation of Kazhdan-Lusztig
polynomials which encodes signature information in addition to composition
factor multiplicities and Jantzen filtration level.  Careful consideration
of Gabber and Joseph's proof of Kazhdan and Lusztig's recursive formula for
computing Kazhdan-Lusztig polynomials and an application of Jantzen's 
determinant formula lead to a recursive formula for the 
signed Kazhdan-Lusztig polynomials.  We use these polynomials
to compute the signature of an invariant Hermitian form on an irreducible
highest weight module.  Such a formula has applications to unitarity testing.
\end{abstract}
\maketitle

\numberwithin{equation}{subsection}
\numberwithin{theorem}{subsection}
\section{Introduction}
\subsection{The Unitary Dual Problem}
In the 1930s, I.M.\ Gelfand introduced a broad programme in abstract 
harmonic analysis which would permit the transfer of difficult problems in 
areas as distinct from analysis as topology to more tractable problems in 
algebra.  Fourier analysis is just one incarnation of this programme.  An 
unresolved component in Gelfand's programme is the classification
of the irreducible unitary representations of a group,
known as the {\it unitary dual problem}.

In the case of a real reductive
Lie group, the problem is equivalent to identifying all irreducible
Harish-Chandra modules which admit a positive definite invariant Hermitian
form.  As Harish-Chandra modules may be constructed via an algebraic
method introduced by Zuckerman in 1978 known as cohomological induction,
it is of interest to study signatures of invariant Hermitian forms
on cohomologically induced modules and to understand how positivity can
fail.

Cohomological induction is a two-step process in which we compose an induction
functor with a Zuckerman functor $\Gamma^i$.  The intermediate module in
cohomological induction is a generalized Verma module which admits an 
invariant Hermitian form if the module to which induction was applied
admits an invariant Hermitian form.  Formulas for the signatures of invariant 
Hermitian forms on these intermediate modules may be used to compute 
signatures of forms on corresponding cohomologically induced modules 
(eg.\ \cite{W}).  This motivates the study of invariant Hermitian forms on 
Verma modules (see \cite{W}, \cite{Y}) and irreducible highest weight modules.

\subsection{Overview}
Let $\g_0$ be a real semisimple Lie algebra, $\theta$ a Cartan involution of 
$\g_0$, $\g_0 = \mfk_0 \oplus \mfp_0$ the corresponding Cartan decomposition, 
and $\h_0 = \mft_0 \oplus \mfa_0$ a $\theta$-stable Cartan subalgebra and 
corresponding Cartan decomposition (recall that every Cartan subalgebra is 
conjugate to one which is $\theta$-stable).  We drop the subscript $0$ to 
denote complexification.  A Hermitian form $\left< \cdot, \cdot \right>$ on a 
$\g$-module $V$ is {\bf invariant} if it satisfies
$$\left< Xv, w \right> + \left< v, \bar{X}w \right> = 0$$
for every $X \in \g$ and every $v, w \in V$, where $\bar{X}$ denotes the 
complex conjugate of $X$ with respect to the real form $\g_0$ of $\g$.   In 
this paper, we develop a formula for the signature character of an invariant 
Hermitian form on an irreducible highest weight module of regular 
infinitesimal character 
when it exists and $\h$ is compact.  In the process, we 
present a survey of results concerning Verma modules and the 
Bernstein-Gelfand-Gelfand category $\cO$.

We describe the organization of this paper.
In section \ref{IrredVerma_section}, we review formulas for signature 
characters of invariant Hermitian forms on Verma modules and their relation 
to the Jantzen filtration.  In section \ref{Verma_section}, we give a brief 
survey of the theory of Verma modules, define signed Kazhdan-Lusztig 
polynomials, and then express the signature character of an invariant 
Hermitian form on an irreducible highest weight module in terms of these 
polynomials (Theorem \ref{hwm_formula}).  In section \ref{SignedKL_section},
we develop recursive formulas for computing signed Kazhdan-Lusztig 
polynomials when $\h$ is compact (Theorem \ref{signedKL_formula}).  We begin 
by treating the elementary cases.  Next, we discuss category $\cO$, Jantzen's 
translation functors, Jantzen's determinant formula, and coherent continuation 
functors; we adapt classical results for contravariant forms to invariant 
Hermitian forms.  Finally, we show how Gabber and Joseph's proof of the 
remaining ``difficult'' recursive formula for computing Kazhdan-Lusztig 
polynomials may be modified to complete a set of recursive formulas which may 
be used to compute signed Kazhdan-Lusztig polynomials.  In section 
\ref{Examples_section}, we compute some examples.  Throughout this paper, 
the results typically only require existence of non-zero invariant Hermitian 
forms (thus $\h$ must be maximally compact in addition to $\theta$-stable).  
However, we only know how to compute signed Kazhdan-Lusztig polynomials when 
$\h$ is compact because in this case, it is easy to derive formulas for 
various quantities from analogous formulas for contravariant forms.  We 
impose the additional condition of compactness on $\h$ in:  subsection 
\ref{EasyRecursion_subsection},
Lemma \ref{form_comparison}, Proposition \ref{form_comp}, and any 
computations of inner products from 
subsection \ref{CoherentContinuation_subsection} to the end of section 
\ref{SignedKL_section}.

\subsection*{Acknowledgements}
I would like to thank David Vogan for many helpful discussions, Joel 
Kamnitzer and Hannah Wachs for their hospitality, Alexander Postnikov 
for asking about irreducible highest weight modules during my thesis 
defence, and the referee for many helpful suggestions.

\section{Forms on Verma Modules and Filtrations}
\label{IrredVerma_section}

\subsection{The Signature of the Shapovalov Form on Irreducible Verma Modules}
We review the contents of \cite{Y}.

Let $\mfb = \h \oplus \mfn$ be a Borel subalgebra of $\g$.
Let $\Delta^+(\g,\h)$ be the set of positive roots determined by 
$\mfb$ and let $\rho$ be one half the sum of the positive roots.  For $\lambda 
\in \h^*$, let $M( \lambda ) = U( \g ) \otimes_{U(\mfb)} \bbC_{\lambda 
- \rho}$ be the Verma module of highest weight $\lambda - \rho$.  We choose 
a generator $v_{\lambda-\rho}$ for $M( \lambda )$.  Given $\mu \in \h^*$, 
$\theta \mu$ is defined by $(\theta\mu)(H) = \mu( \theta^{-1}
H )$ for every $H$ in $\h^*$.  The complex conjugate of $\mu$, 
$\bar{\mu}$, is defined by $\bar{\mu}(H) = \overline{\mu(\bar{H})}$.
We have the relation $\theta \mu  = - \bar{\mu}$ for $\mu \in \Lambda_r$.
The Verma module $M( \lambda )$ admits an invariant Hermitian form $\left< 
\cdot, \cdot \right>_\lambda$, which is 
unique up to a real scalar, when $\h$ is maximally compact, 
$\theta (\Delta^+( \g, \h)) = \Delta^+(\g,\h)$ (recall that $\h$ is 
$\theta$-stable), and $\lambda$ takes imaginary 
values on $\h_0$.  Normalized so that $\left< v_{\lambda-\rho}, v_{\lambda 
-\rho} \right>_\lambda = 1$, it is known as the Shapovalov form.  Henceforth, 
we shall take $\h$ and $\mfb$ to have the properties required for existence 
of invariant Hermitian forms.

Define $X^* = -\bar{X}$ for $X \in \g$, and extend $\cdot^*$ to an 
anti-involution of $U(\g)$ via $(xy)^* = y^*x^*$ for $x, y \in U( \g )$.  
Then
\begin{equation}
\label{hermitian_formula}
\left< xv_{\lambda-\rho}, yv_{\lambda-\rho} \right>_\lambda = (\lambda-\rho)
( p( y^*x ) ) 
\end{equation}
where $p : U( \g ) \to U( \h )$ is defined to be projection under the direct 
sum $U(\g) = U( \h) \oplus ( \mfn^{op}U( \g) + U( \g )\mfn )$. 

Due to invariance, the Shapovalov 
form pairs the $\lambda - \mu - \rho$ weight space of $M( \lambda )$ with 
the $\lambda - \theta\mu - \rho$ weight space.
Recalling that $\theta \mu= -\bar{\mu}$ for $\mu \in \Lambda_r$, 
we see that these are two distinct 
weight spaces when $\mu$ is not imaginary, and they are the same weight 
space if $\mu$ is imaginary.  Finite dimensionality of these spaces allows us 
to discuss determinants and signatures of the Shapovalov form.

A modification of the classical (invariant bilinear) Shapovalov determinant 
formula shows that when $\mu$ is imaginary so that the weight space 
$M( \lambda)_{\lambda-\mu-\rho}$ is paired with itself, the determinant of a 
matrix representing the Shapovalov form on $M( \lambda )_{\lambda-\mu-\rho}$ 
is
$$ \prod_{\alpha \in \Delta^+(\g,\h)} \prod_{n=1}^{\infty} \left( \left( 
\lambda, \alpha^\vee \right) - n \right)^{P(\mu - n\alpha)}$$
up to multiplication by a scalar determined by the basis chosen.  $P$ denotes 
Kostant's partition function.

When $\mu$ is not imaginary (i.e. $\mu \neq \theta\mu$), the determinant of a 
matrix representing 
$\left< \cdot, \cdot \right>_\lambda$ on $M( \lambda )_{\lambda - \mu - \rho}
\oplus M( \lambda )_{\lambda - \theta\mu -\rho}$ is
$$ \prod_{\alpha \in
\Delta^+(\g,\h)} \prod_{n=1}^{\infty} \left( \left( \lambda, \alpha^\vee
\right) - n \right)^{P(\mu - n\alpha)}\left( \left( \lambda, \alpha^\vee
\right) - n \right)^{P(\theta\mu - n\alpha)}$$
up to multiplication by a scalar.  

Since $\left< M( \lambda )_{\lambda-\nu-\rho}, M( \lambda )_{\lambda - \nu
-\rho} \right>_\lambda = 0$ when $\nu$ is not imaginary, therefore for
non-imaginary $\mu$, 
the number of positive eigenvalues of a matrix representing the 
Shapovalov form on $M( \lambda )_{\lambda-\mu-\rho} \oplus M( \lambda )_{
\lambda -\theta\mu -\rho}$ equals the number of negative eigenvalues of that 
matrix (see Sublemma 3.18 of \cite{V3} or p.7 of \cite{Y}).  Thus the 
signature of the Shapovalov form on 
$M( \lambda )$ may be 
recorded in a formal sum called the {\bf signature character} as follows:
$$ch_s M( \lambda ) = \sum_{ \stackrel{\mu \in \Lambda_r^+}{\mu 
\text{ imaginary}}} ( p( \mu ) - q( \mu ))e^{\lambda-\mu-\rho}$$
where the signature of the form on $M( \lambda )_{\lambda-\mu-\rho}$ is 
$( p(\mu), q( \mu ) )$.  Note that $p( \mu ) + q( \mu )$ is the dimension of 
$M( \lambda )_{\lambda - \mu - \rho}$, whence the usual character formula 
when all the roots are imaginary (i.e. $\h$ is compact) is

$$ ch \, M( \lambda ) = \sum_{\mu \in \Lambda_r^+} ( p( \mu ) + q( \mu ) )
e^{\lambda-\mu-\rho}.$$

The radical of the Shapovalov 
form is the unique maximal submodule of $M( \lambda )$, whence Verma modules 
are reducible precisely when the Shapovalov form is degenerate.  
Thus Verma modules $M(\lambda)$ are 
reducible precisely on the affine hyperplanes $H_{\alpha,n} := \{ \lambda \, 
| \, (\lambda, \alpha^\vee) = n \}$ where $\alpha$ is a positive root and $n 
\in \bbZ^+$.  Within any region defined by these reducibility 
hyperplanes, the signature of the Shapovalov form cannot change because the 
form remains non-degenerate.

The largest region for which the signature does not change is the
intersection of the negative open half spaces
$$\left( \bigcap_{\alpha \in \Pi}H_{\alpha,1}^-  \right)
\bigcap H_{\widetilde{\alpha}, 1}^-$$
with $i\h_0^*$,
where $H_{\alpha,n}^- := \{ \lambda \, | \, (\lambda, \alpha^\vee) < n \}$,
$\widetilde{\alpha}^\vee$ is the highest coroot, and $\Pi$ is the set of
simple roots corresponding to our choice of $\Delta^+$.  In \cite{W}, Wallach 
used an asymptotic argument to calculate the signature character of Verma 
modules $M( \lambda )$ with $\lambda$ in this region,
which we refer to as the Wallach region:
\begin{theorem}
\label{Wallach_formula}
Let imaginary $\lambda$ satisfying $\lambda |_{\mfa_0} \equiv 0$
be in the Wallach region.  Then the signature character of 
the Shapovalov form $\left< \cdot, \cdot \right>_\lambda$ on $M(\lambda)$ is
\begin{equation*}
ch_s M(\lambda ) = \frac{e^{\lambda-\rho}}{\displaystyle{\prod_{ \alpha \in
\Delta^+( \mfp, \mft ) } \left( 1 - e^{-\alpha} \right)
\prod_{ \alpha \in \Delta^+( \mfk, \mft )
} \left( 1 + e^{-\alpha} \right)} }.
\end{equation*}
\end{theorem}

(In section 2 of \cite{Y}, we discuss compatibility with the definition of 
the signature character above.)
The form taken by reducibility hyperplanes suggests searching for a formula 
for the signature character in other regions which uses the affine Weyl 
group.  In \cite{Y}, 
we defined 
$$R( \lambda ) :=
\sum_{\mu \in \Lambda_r^+} {c_\mu e^{\lambda - \mu}}$$  
for some constants $c_\mu$ given by Wallach's formula 
so that $R(\lambda )$ is the signature character of the Shapovalov form 
$\left< \cdot, \cdot \right>_\lambda$ when $\lambda$ lies in the Wallach 
region.  We also defined and computed
$$R^A( \lambda ) := \displaystyle{ \sum_{\mu \in \Lambda_r^+} c^A_\mu
e^{\lambda - \mu}}$$ 
so that $R^A( \lambda )$ is the signature character of the Shapovalov form
$\left< \cdot, \cdot \right>_\lambda$ when
$\lambda$ belongs to the alcove $A$.  We chose the fundamental Weyl 
chamber $\mathfrak{C}_0$ to be the antidominant chamber and the fundamental 
alcove $A_0$ to be in the antidominant Weyl chamber.
We computed the signature of the Shapovalov form within an alcove $A$ by 
taking a path from $A$ to a specific alcove in the Wallach region and 
computing changes to the signature character for each hyperplane crossed.  
We refer the reader to Theorem \ref{Verma_formula} of this paper or Theorem 
4.6 of \cite{Y} for a partial statement of 
the formula and to Theorem 6.12 of \cite{Y} for the full formula.

\subsection{The Jantzen Filtration}
As it is fundamental to our study of irreducible highest weight modules, we
review the main tool used in computing how signatures change across a 
hyperplane: the Jantzen filtration.  The Jantzen filtration corresponding to 
an analytic family of Hermitian forms $\left< \cdot, \cdot \right>_t$, for 
$t \in (- \delta, \delta)$, on a finite dimensional vector space $E$ 
is the sequence of subspaces
$$E = E^{\left<0\right>} \supset E^{\left<1\right>} \supset \cdots 
\supset E^{\left<N\right>} = \lbrace 0 \rbrace$$
where $e \in E^{\left<n\right>}$ for $n \geq 0$ if there exists an analytic 
function $f_e :
( -\varepsilon, \varepsilon ) \to E$ for some $\varepsilon > 0$ such that
\begin{enumerate}
\item $f_e( 0 ) = e$
\item $\left< f_e(t), e' \right>_t$ vanishes to order at least $n$ at $t=0$
for any $e' \in E$.
\end{enumerate}
As weight spaces of a Verma module are finite dimensional, by taking an 
analytic path $\lambda_t: (-\delta, \delta) \to i\h_0^*$ and corresponding 
Hermitian forms $\left< \cdot, \cdot \right>_{\lambda_t}$, we may discuss 
Jantzen filtrations of Verma modules.

For $e, e' \in E^{\left<n\right>}$, define
$$ \left< e, e' \right>^n = \lim_{t \to 0} \frac{1}{t^n} \left< f_e( t ),
f_{e'}(t) \right>_t$$
which is independent of choice of $f_e$ and $f_{e'}$.  Then
\begin{theorem} (\cite{V3}, Proposition 3.3)
\label{jantzen_filtration}
The form $\left< \cdot, \cdot \right>^n$ on $E^n$ is Hermitian with radical
$E^{\left<n+1\right>}$, and therefore it induces a non-degenerate Hermitian 
form on $E_{\left<n\right>} := E^{\left<n\right>} / E^{\left<n+1\right>}$, 
which we also denote $\left< \cdot, \cdot \right>^n$.
Let $(p_n, q_n)$ be the signature of $\left< \cdot, \cdot \right>^n$,
$(p,q)$ be the signature of $\left< \cdot , \cdot \right>_t$ for $t \in (0,
\delta)$, and $(p',q')$ be the signature of $\left< \cdot , \cdot
\right>_t$ for $t \in (-\delta,0)$.
Then
\begin{eqnarray*}
( p, q)  &=& \left( \sum_n p_n , \sum_n q_n  \right) \qquad \text{and}\\
( p', q')  &=& \left( \sum_{n \text{ even}} p_n + \sum_{n \text{ odd}} q_n, 
 \sum_{n \text{ odd}} p_n + \sum_{n \text{ even}} q_n \right).
\end{eqnarray*}
\end{theorem}
Consider an analytic path $\lambda_t$ such that $\lambda_0$ lies in exactly 
one reducibility hyperplane, $H_{\alpha,n}$, and $\lambda_t$ does not lie in 
any reducibility hyperplane for $t \neq 0$.  Then $M(\lambda_0 )$ 
has a unique proper non-trivial submodule:  $M( \lambda_0 - n\alpha) = 
M( s_\alpha \lambda_0 )$.  It 
lies in an odd level of the Jantzen filtration, and therefore as one crosses 
the hyperplane $H_{\alpha,n}$, the signature changes by the signature of 
an invariant Hermitian form on $M(\lambda_0 -n\alpha)$ and thus by plus or 
minus the signature of $\left< \cdot, \cdot \right>_{\lambda_0 -n\alpha}$.
We write this as
$$R^{A}( \lambda ) = R^{A'}( \lambda ) + 2 \varepsilon( A,
A') R^{A - n\alpha}( \lambda - n \alpha )$$
for adjacent alcoves $A$ and $A'$ separated by the hyperplane 
$H_{\alpha,n}$.  $\varepsilon( A, A' )$ is a function of the Weyl 
chamber containing $A$ and $A'$, $\alpha$, and $n$.  Its value may be found 
in \cite{Y}.  The relation above leads to the inductive formula for the 
signature character of the Shapovalov form on an irreducible Verma module 
which was mentioned previously.

Results which we wish to use are formulated in terms of the canonical 
Jantzen filtration and so we review this classical concept and investigate 
relations to our version of the Jantzen filtration.

We use the setup of \cite{GJ} since we will follow sections of it 
closely.  Let $\lbrace X_\alpha, Y_\alpha \, | \, \alpha \in \Delta^+( \g, \h
) \rbrace \cup \lbrace H_\alpha \, | \, \alpha \in \Pi \rbrace$ be the 
Chevalley basis for $\g$.  We let $\g_\bbZ$ be the integer span of the 
Chevalley basis.  It is a Lie algebra.  $\h_\bbZ$, $\mfn_\bbZ$, 
$\mfn^{op}_\bbZ$, and $\mfb_\bbZ$ are the obvious analogues.  We let $A$ be 
the local ring $\bbC\left[t\right]_{(t)}$.  For the Lie algebra $a_\bbZ$, 
define $a_A$ to be $a_\bbZ \otimes_{\bbZ}A$.  Jantzen defines $U(\mfn^{op}
)_\bbZ$ to be the $\bbZ$-subalgebra generated by $Y_{\alpha,n} = Y_\alpha^n
/n!$ where $\alpha \in \Delta^+( \g,\h)$ and $n \geq 0$.  For a highest weight 
module $E$ with primitive generator $v$ for which $\left< v, v \right> = 1$, 
define $E_\bbZ$ to be $U(\mfn^{op})_\bbZ v$.

Given $\lambda \in \h_A^*$, we let $A_\lambda$ be the one dimensional 
$U( \h_A )$-module on which $H \in \h$ acts by multiplication by $\lambda(H)$.
Extending $A_\lambda$ to a $U( \mfb_A)$-module by allowing $X \in \mfn_A$ to 
act by zero, we define the Verma module over $U( \g_A )$ by 
$M( \lambda )_A := U( \g_A ) \otimes_{U( \mfb_A)} A_{\lambda - \rho}$.

There is an involutive antiautomorphism $\sigma$ of $\g$ so that $\sigma( 
H ) = H$ for all $H \in \h$ and such that $\sigma( X_\alpha ) = Y_\alpha$ 
for every positive root $\alpha$.  It may be extended to an antiautomorphism 
of $U( \g )$ in the same way that $\cdot^*$ was.  This leads to the 
canonical contravariant 
forms, denoted by $( \cdot , \cdot )$, on $U(\g)$ and $U( \g_A)$ Verma modules 
with the defining properties that the forms are symmetric, bilinear, and 
$( xv, w ) = ( v, \sigma( x )w )$ for all $x$ in the universal enveloping 
algebra and all elements $v, w$ of the Verma module.  Thus $( \cdot, \cdot )$ 
on $M( \lambda )$ or $M( \lambda )_A$ satisfies
\begin{equation}
\label{contravariant_formula}
( xv_{\lambda-\rho}, yv_{\lambda-\rho}) = (\lambda-\rho)( p( \sigma(y)x ) ).
\end{equation}
Compare this with (\ref{hermitian_formula}).

Consider the $U( \g_A)$-module $M = M(\lambda + \delta t)_A$ where $\lambda \in
\h^*$ and $\delta \in \h^*$ are regular and imaginary.  The Jantzen filtration 
of $M( \lambda + \delta t)_A$ is defined to be
$$ M^{(0)} \supset M^{(1)} \supset \cdots \supset M^{(N)} = \lbrace 0 \rbrace$$
where $M^{(j)} = \lbrace v \in M \, |
\, (v,w) \in (t^j) \, \forall \, w \in M \rbrace$.  
This is a filtration of $M( \lambda + \delta t)_A$ by 
$U( \g_A )$ modules.  We get a filtration on the module $\bar{M} := 
M /{tM}$, which is isomorphic to the $U( \g )$ Verma module $M( \lambda )$, 
via $\bar{M}^{(j)} = M^{(j)} / (tM \cap M^{(j)})$.  It is the usual 
Jantzen filtration of $M( \lambda )$ and {\it does not depend on the value of 
regular $\delta$}:  in \cite{BA}, Barbasch showed for an arbitrary 
non-degenerate deformation direction that the Jantzen filtration 
coincides with the socle filtration.
We define $M_{(j)}$ and $\bar{M}_{(j)}$ as in the Hermitian case.

\begin{lemma}
\label{Jantzen_same}
Let $M(\lambda)^{\left<j\right>}$ be the $j^{\text{th}}$ level of the 
Jantzen filtration defined by the path 
$\lambda_t = \lambda + \delta t$.  Then
$$M(\lambda)^{(j)} = M(\lambda)^{\left<j\right>}.$$
\end{lemma}
\begin{proof}
Note that 
$(\lambda + \delta t - \rho) (y^*x)
=(\lambda + \delta t - \rho) (\sigma(\sigma(y^*))x) $
so that
$$\left< xv_{\lambda +  \delta t -\rho}, y v_{
\lambda + \delta t - \rho} \right>_{\lambda + \delta t} 
= (x v_{\lambda + \delta t
- \rho}, \sigma(y^*)v_{\lambda + \delta t - \rho}).$$
The lemma now follows from the two definitions of the Jantzen filtration and 
the observation that $\cdot^*$ and $\sigma$ are bijections from $U( \g )$ 
to $U( \g )$.
\end{proof}

Henceforth, bar will denote specialization at $t=0$.
For the remainder of this paper, we use the classical Jantzen filtration of a 
Verma module and use $j$ interchangeably with $(j)$ and with $\left< j \right>$.

\section{Verma modules and Kazhdan-Lusztig polynomials}
\label{Verma_section}
\subsection{A Brief Overview of Verma Modules}
The structure of Verma modules has been studied by a number of people 
(eg.  \cite{VE}, \cite{BGG2}, \cite{DL}, \cite{J}, \cite{BB}). 
\begin{theorem} (cf. \cite{D} Theorem 7.6.6)
$$\dim \Hom_\g( M( \lambda ), M( \mu ) ) \leq 1 \text{ for all }
\lambda, \mu \in \h^*.$$
\end{theorem}
\begin{theorem} (Bernstein-Gelfand-Gelfand, \cite{D} Theorem 7.6.23) 
For $\lambda, \mu
\in \h^*$,
\begin{eqnarray*}
M( \mu ) \subset M( \lambda) \iff &
\exists \,\alpha_1, \cdots, \alpha_m \in \Delta^+(\g, \h)
\text{ such that } \\
& \lambda \geq s_{\alpha_1}\lambda \geq \cdots \geq
s_{\alpha_m} \cdots s_{\alpha_1} \lambda = \mu.
\end{eqnarray*}
(Recall that for $\mu_1, \mu_2 \in \h^*, \mu_1 \leq \mu_2$ if and only if 
$\mu_2 - \mu_1 \in \Lambda_r^+$.)
\label{BGG}
\end{theorem}
\begin{remark}
The above conditions may not be equivalent to
$\mu \in W \lambda \text{ and } \mu \leq \lambda$.
\end{remark}
Verma modules have finite composition series.  The composition factors of 
$M(\lambda)$ are $L(\mu)$ where $M( \mu ) \subset M( \lambda )$ (cf. 
\cite{D}, Theorem 7.6.23).
In \cite{DL}, Deodhar and Lepowsky showed that although $\dim \Hom_\g( M( 
\mu ), M( \lambda ) ) \leq 1$, it is possible for a composition factor of a 
Verma module to have multiplicity greater than one.  In \cite{KL}, Kazhdan 
and Lusztig defined polynomials $P_{x,y}$ for $x,y \in W$ known as 
Kazhdan-Lusztig polynomials.  They famously conjectured that for $\lambda$ 
antidominant and regular and for $x$ and $y$ in the integral Weyl group 
$W_\lambda$ with longest element $w_\lambda$, the polynomials give the 
multiplicity of $L( y\lambda)$ as a composition factor of $M( x\lambda )$:
\begin{equation*}
\left[ M( x \lambda ): L( y \lambda ) \right] = P_{w_\lambda x, 
w_\lambda y}(1)
\end{equation*}
from which we obtain
\begin{equation}
\label{KL_mult}
ch \, M( x\lambda ) = \sum_{y \leq x} P_{w_\lambda x, w_\lambda y}(1) 
ch \, L( y \lambda ).
\end{equation}
Furthermore, 
the multiplicity of $L( y\lambda )$ in the $j^{\text{th}}$ level of the 
Jantzen filtration of $M( x\lambda )$ is encoded by Kazhdan-Lusztig 
polynomials:
\begin{equation}
\left[M( x\lambda )_j: L(y\lambda ) \right] = \text{the coefficient of } 
q^{(\ell(x) - \ell(y) - j )/2} \text{ in } P_{w_\lambda x, w_\lambda y}.
\end{equation}

A proof of the Kazhdan-Lusztig Conjecture was perhaps the most important 
open problem in representation theory in the early 1980s.
In \cite{V4}, Vogan showed that semisimplicity  of $U_\alpha L(x \lambda)$, 
where $U_\alpha L(x \lambda )$ is defined to be the cohomology
of the complex $ 0 \rightarrow L( x \lambda ) \rightarrow \theta_\alpha
L( x \lambda ) \rightarrow L( x \lambda ) \rightarrow 0$,
implies the Kazhdan-Lusztig Conjecture.  In \cite{GJ}, Gabber and Joseph
proved that Vogan's Conjecture follows from Jantzen's Conjecture: 
$$ M( x \lambda )^j = M( xs_\alpha \lambda )^{j+1} \cap M( x \lambda )$$
for $j \geq 0$, $x \in W_\lambda$, $xs_\alpha > x$, and $\alpha \in \Pi$ such 
that $( \lambda, \alpha^\vee ) \in \bbZ$.
Brylinski, Kashiwara, Beilinson, and Bernstein were able to prove the 
Kazhdan-Lusztig Conjecture by studying the relation between Kazhdan-Lusztig 
polynomials and Deligne, Goresky, and MacPherson's intersection 
cohomology (\cite{KB}, \cite{BB2}).  Beilinson and Bernstein subsequently 
proved Jantzen's Conjecture (\cite{BB}) using stronger versions of these 
techniques.

\subsection{A formula for $ch_s L( x\lambda )$ in terms of signed 
Kazhdan-Lusztig polynomials}
Because the radical of the Shapovalov form  on 
$M( x \lambda )$ is $M( x \lambda )^1$ and $L( x \lambda ) = M( x\lambda )_0
= M(x\lambda)/M(x\lambda)^1$, 
therefore the Shapovalov form on the Verma module descends to an invariant 
Hermitian form, which we also call the Shapovalov form, on the irreducible 
highest weight module $L( x\lambda )$.  Their signatures differ only by zero 
eigenvalues.  We write $ch_s L( x\lambda )$ for the signature character of 
the Shapovalov form on $L( x\lambda )$.  (We implicitly assume here that 
$x\lambda$ is imaginary.)

From equation (\ref{KL_mult}), one obtains the inversion formula
\begin{equation*}
ch \, L( x \lambda ) = \sum_{y \leq x} (-1)^{\ell( x ) - \ell( y) }
P_{y, x}( 1 )\, ch \, M( y \lambda )
\end{equation*}
(cf. \cite{KL}).  
Because we do not know the value of $ch_s M(y \lambda )$ when $M( y \lambda )$ 
is reducible, we cannot compute signature characters from the above formula.  
However, we may make use of our knowledge of signature characters for 
alcoves which contain $x\lambda$ in their closures.  
To illustrate this, consider the 
simple example of $x\lambda$ such that only the adjacent alcoves $A$ and $A'$ 
contain $x\lambda$ in their closures.  Let $H_{\alpha,n}$ be the reducibility
hyperplane containing $x\lambda$.  Recall that as one crosses the hyperplane 
$H_{\alpha,n}$ at $x\lambda$, the signature changes by the signature of 
$M( x\lambda 
-n\alpha) = M( s_\alpha x\lambda ) = L( s_\alpha x \lambda )$.
We conclude that the signature characters for the 
alcoves $A$ and $A'$ evaluated at the point $x\lambda$ are $ch_s L( x\lambda ) 
\pm ch_s L(s_\alpha x\lambda )$ in some order so that $\frac{1}{2}
\left( R^A( x\lambda ) + 
R^{A'}(x\lambda )\right) = ch_s L( x\lambda )$.  We have formulas for $R^A$ 
and for $R^{A'}$, and so we have expressed $ch_s L( x\lambda )$ in terms of 
known quantities.

We now consider the general case.  Take the path $\lambda_t = x\lambda + 
\delta t$, where $\delta$ is regular and imaginary, and consider the 
Jantzen filtration of $M( x\lambda )$ which it defines.
Now $M( x\lambda )_j$ is semisimple (cf. \cite{GJ} Theorem
4.8 (ii) ) and  $\left< \cdot, \cdot \right>_j$ is a non-degenerate invariant 
Hermitian form on $M( x\lambda )_j$.  The contribution to the signature 
character 
of $\left< \cdot, \cdot \right>_j$ by a particular irreducible constituent 
$L( y\lambda )$ of $M( x\lambda )_j$ is either the signature character of 
the Shapovalov form, the negative of it, or zero because $L( y\lambda )$ is 
paired with $L( \theta y\lambda)$ (which may be another copy of $L( 
y\lambda )$).  Recording 
which of the three choices occurs for each composition factor with $+1$, $-1$, 
or $0$, we have
\begin{equation}
ch_s \left< \cdot, \cdot \right>_j = \sum_{y \leq x} a^{\lambda, 
\delta}_{ w_\lambda x, w_\lambda y, j}
ch_s L( y\lambda )
\end{equation}
for some integers $a^{\lambda, \delta}_{ w_\lambda x, w_\lambda y, j}$.
Since signatures cannot change in the interior of an alcove, we will let 
$w( \delta ) \in W_\lambda$ be such that $\delta \in w(\delta) \mathfrak{C}_0$ 
and we will write $a^{\lambda,w(\delta)}_{w_\lambda x, w_\lambda y, j}$
in place of $a^{\lambda,\delta}_{w_\lambda x, w_\lambda y, j}$  We record 
these integers in polynomials
$$ P^{\lambda,w}_{w_\lambda x, w_\lambda y}(q) := \sum_{j \geq 0} a^{\lambda,w}
_{w_\lambda x, w_\lambda y, j} q^{\frac{\ell(x)-\ell(y)-j}{2}}$$
which we call the {\bf signed Kazhdan-Lusztig polynomials}.  We remind the 
reader that the signed Kazhdan-Lusztig polynomials above are indexed by a 
regular antidomininant weight $\lambda$ and by elements $w$, $x$, and $y$ of 
$W_\lambda$.

\begin{remark}
Note that the usual Kazhdan-Lusztig polynomials are defined in the same way, 
but with contributions of $+1$ to coefficients for every composition factor 
rather than $+1$, $-1$, or $0$.  Therefore
$$\left| a^{\lambda,w}_{w_\lambda x, w_\lambda y, j} \right| \leq 
\left[ M( w\lambda)_j : L( y \lambda )\right] = \text{coefficient of }
q^{(\ell(x)-\ell(y)-j)/2} \text{ in } P_{w_\lambda x, w_\lambda y}.$$
\end{remark}

Let $A( x \lambda, w(\delta))$ be the alcove containing $x\lambda + \delta t$ 
for regular, imaginary $\delta$ and for small $t > 0$.  Using our formulas 
above and Theorem \ref{jantzen_filtration}:
$$R^{A(x\lambda,w(\delta))}(x\lambda) = \sum_j ch_s \left< \cdot, \cdot 
\right>_j = \sum_{y \leq x} P^{\lambda,w(\delta)}_{w_\lambda x, w_\lambda y}
(1) ch_s L( y\lambda ).$$

\begin{remark}
When $y\lambda$ is not imaginary, neither $M( y\lambda )$ nor $L(y\lambda)$
admit non-zero invariant Hermitian forms, whence $ch_s L(y\lambda)$ is 
undefined.  However, in this case, each $L( y\lambda)$ is paired with some 
$L( \theta y\lambda )$ , giving us $a^{\lambda,w(\delta)}_{w_\lambda x,
w_\lambda y, j} = 0$ and $P^{\lambda, w(\delta)}_{w_\lambda x, w_\lambda y}
= 0$ (see our discussion in the previous section concerning non-imaginary 
weights).  The sum above ought to be over $y$ such that $y\lambda$ is 
imaginary to avoid abuse of notation.
\end{remark}
Observe that $P^{\lambda,w(\delta)}_{w_\lambda x, w_\lambda x}(q) = 
a^{\lambda,w(\delta)}_{w_\lambda x, w_\lambda x,0} = 1$, and so
$$ ch_s L( x\lambda ) = R^{A(x\lambda ,w(\delta))}(x\lambda) - \sum_{\stackrel
{y < x}{y\lambda \text{ imaginary}}}
P^{\lambda,w(\delta)}_{w_\lambda x, w_\lambda y}(1) ch_s L( y\lambda ).$$

Applying this formula recursively, we arrive at:
\begin{theorem}
\label{hwm_formula}
If $\lambda \in \h^*$ is regular and antidominant, then for $x \in W_\lambda$ 
such that $x\lambda$ is imaginary and for any $w \in W_\lambda$:
$$ch_s L( x\lambda ) = \sum_{\stackrel{y_1 < \cdots < y_j=x} 
{y_k\lambda \text{'s imaginary}}}
(-1)^{j-1} \left( \prod_{i=2}^{j} P^{\lambda,w}_{w_\lambda y_i, w_\lambda 
y_{i-1}}(1) \right)
R^{A(y_1\lambda , w)}(y_1 \lambda).$$
\end{theorem}

We recall the formula for $R^{A(y_1\lambda, w)}(y_1 \lambda)$:
\begin{theorem}
\label{Verma_formula}
(\cite{Y}, Theorems 4.6 and 6.12)  Let $\Delta_i^+( \g, \h)$ be the set of 
imaginary roots in $\Delta^+( \g, \h)$.  Subscripts or superscripts $i$ will 
refer to objects associated with $\Delta_i^+( \g, \h)$.  We will assume that 
everything (simple roots, reducibility hyperplanes, etc.) in this theorem 
is associated to the root system of imaginary roots.  Choose the fundamental 
alcove $A_0^i$ of $W_a^i$ and the fundamental chamber $\mathfrak{C}_0^i$ of 
$W_i$ to contain $-\rho_i$.  Let $\bar{\cdot}: W_a^i \to W_i$ be the 
homomorphism arising from the semidirect product structure $W_a^i = W_i \ltimes
\Lambda_i$.  Given $a \in W_a^i$, let $\widetilde{a} \in W_i$ be such that 
$aA_0^i \in \widetilde{a}\mathfrak{C}_0^i$.  Let $aA_0^i = C_0 
\stackrel{r_1}{\to} C_1 \stackrel{r_2}{\to} \cdots \stackrel{r_\ell}{\to}
C_\ell = \widetilde{a}A_0^i$ be a path from $aA_0^i$ to $\widetilde{a}A_0^i$.
Then for imaginary $\lambda \in A_0^i$:
\begin{eqnarray*}
ch_s M( \lambda ) |_{\mfa_0} &=& \lambda |_{\mfa_0} \qquad \text{and} \\
ch_s M( \lambda ) |_{\mft_0} &=& R^{aA_0^i}( \lambda |\mft_0 ) \\
&=& \sum_{\stackrel{S=\lbrace i_1 < \cdots i_k \rbrace}
{\subset \lbrace 1, \ldots, \ell \rbrace}} \varepsilon( S ) 2^{|S|}
\frac{e^{\overline{r_{i_1}r_{i_2}}\cdots \overline{r_{i_k}}r_{i_k} r_{i_{k-1}}
\cdots r_{i_1} \lambda |_{\mft_0} - \rho}}
{ \prod_{\alpha \in \Delta^+( \mfp, \mft )}( 1 - e^{-\alpha} )
\prod_{\alpha \in \Delta^+( \mfk, \mft )}( 1 + e^{-\alpha} )}
\end{eqnarray*}
where
$\varepsilon(S) = \varepsilon( C_{i_1 - 1}, C_{i_1} ) \varepsilon( 
\overline{r_{i_1}}C_{i_2-1}, \overline{r_{i_1}}C_{i_2} ) \cdots 
\varepsilon( \overline{r_{i_1}} \cdots \overline{r_{i_{k-1}}} C_{i_k-1},
\overline{r_{i_1}} \cdots \overline{r_{i_{k-1}}}C_{i_k} )$, $\varepsilon( 
\emptyset ) = 1$, and the formula for $\varepsilon( C, C' )$ for alcoves $C$, 
$C'$ may be found in Theorem 6.12 of \cite{Y}.
\end{theorem}

Since we have formulas for $R^A$ for any alcove $A$,
therefore we can compute  $ch_s L( x\lambda )$ as long as we can 
compute the integers $P^{\lambda,w}_{w_\lambda x, w_\lambda y}(1)$.

\section{Recursive Formulas for computing signed Kazhdan-Lusztig polynomials}
\label{SignedKL_section}
\subsection{Recursive formulas for the easy cases}
\label{EasyRecursion_subsection}
The usual Kazhdan-Lusztig polynomials may be computed via $P_{x,x} = 1$, 
$P_{x,y} = 0$ when $x > y$, and by the recursive formulas:

\begin{itemize} 
\item[a)]
$P_{w_\lambda x, w_\lambda y} = P_{w_\lambda xs, w_\lambda y}$
if $ys > y$ and $x, xs \geq y$, $s$ simple.
\item[a')]
$P_{w_\lambda x, w_\lambda y} = P_{w_\lambda sx, w_\lambda y}$
if $sy > y$ and $x, sx \geq y$, $s$ simple.
\item[b)]
If $y > ys$ then
$$
q^c P_{w_\lambda xs, w_\lambda y} + q^{1-c}P_{w_\lambda x,
w_\lambda y} = 
\begin{array}{c}
\displaystyle{
\sum_{ z \in W_\lambda | zs > z }
{ \mu( w_\lambda z, w_\lambda y) q^{\frac{\ell(z) - \ell(y) + 1}{2}} 
P_{w_\lambda x, w_\lambda z} }} \\ 
+ P_{w_\lambda x, w_\lambda ys}
\end{array}
$$
where $c=1$ if $xs < x$, $c=0$ if $xs>x$, 
and $\mu( w_\lambda z, w_\lambda y)$ is the multiplicity of $L( y \lambda )$
in $M( z \lambda )_1$.
\end{itemize}

The initial conditions for signed Kazhdan-Lusztig polynomials are identical:
$P^{\lambda,w}_{x,x} = 1$ and
$P^{\lambda,w}_{x,y} = 0$ when $x > y$.
We will find the recursive formulas to be different.

We discuss the signed Kazhdan-Lusztig polynomial analogue of case a').  
{\it We assume the Cartan subalgebra $\h$ to be not only 
maximally compact but compact for this subsection}.  
Thus $y\lambda = \theta (y \lambda)$ for 
all $y \in W_\lambda$.  Choose $x,y \in W_\lambda$ and $s=s_\alpha$ simple 
so that $sy > y$ and $sx > x$.
Recall that $\lambda$ is regular antidominant.
We consider the Jantzen filtration 
corresponding to the path $\lambda_t = x\lambda + \delta t$ where $\delta \in 
w\mathfrak{C}_0$ and we restrict our attention to $L( y\lambda )$ in 
the $j+1^\text{st}$ level of $M( sx\lambda )$.  Since $sx > x$, therefore 
$(x\lambda, \alpha^\vee) > 0$, whence $M( x\lambda )$ is a submodule of 
$M( sx\lambda )$.  The Jantzen Conjecture and our recursive formula a') tell 
us that all copies of $L( y \lambda )$ in the $j+1^\text{st}$ level of the 
filtration arise from the submodule $M( x\lambda )$ of $M( sx\lambda )$.  
Similarly, if $xs > x$ and $ys > y$ then all copies of $L( y\lambda )$ in 
the $j+1^{\text{st}}$ level of the Jantzen filtration of $M( xs\lambda )$ 
arise from the submodule $M( x\lambda )$ of $M( sx\lambda )$.
We obtain:
\begin{proposition}  Let $s = s_\alpha$ be a simple reflection and let 
$y < x$.
\begin{itemize}
  \item[a)] $a^{\lambda, w}_{w_\lambda x, w_\lambda y,j+1} = \sgn( -w\rho, 
x\alpha ) \varepsilon( H_{x\alpha,(xs\lambda,x\alpha^\vee)}, xs ) 
a^{\lambda, w}_{w_\lambda x, w_\lambda y,j} $ if 
  $xs > x$ 
  \item[a')] $a^{\lambda, w}_{w_\lambda sx, w_\lambda y,j+1} = \sgn( -w\rho, 
\alpha ) \varepsilon( H_{\alpha,(sx\lambda,\alpha^\vee)}, sx ) a^{\lambda, 
w}_{w_\lambda x, w_\lambda y,j} $ if $ sx > x$.
\end{itemize}
\end{proposition}
\begin{proof}
Recall that $\varepsilon( A, A' )$ is a function of the Weyl chamber 
containing $A$ and $A'$ and the hyperplane which separates them.  We therefore
defined $\varepsilon( H_{\gamma,N}, z )$ for $z \in W$ in \cite{Y}.   Take 
an analytic path $\lambda_t : (-\varepsilon, \varepsilon ) \to i\h_0^*$ in 
the Weyl chamber $z \mathfrak{C}_0$ so 
that $\lambda_t \in H_{\gamma,N}^+$ for $t>0$, $\lambda_t \in H_{\gamma,N}^-$ 
for $t<0$, $H_{\gamma,N}$ is the only reducibility hyperplane containing 
$\lambda_0$, and $M( \lambda_t )$ is irreducible 
for $t \neq 0$.  Let $t_1 \in (0, \varepsilon)$ and let $t_{2} \in 
(-\varepsilon, 0)$.  Recall that 
$$ ch_s M( \lambda_{t_1} ) = e^{\lambda_{t_1} - \lambda_{t_2}} 
ch_s M( \lambda_{t_2} ) + 2 \varepsilon( H_{\gamma,N}, z ) e^{\lambda_{t_1}-
\lambda_0} ch_s M( \lambda_0 -N\gamma )$$
which reflects the change of the signature character by the signature 
character of the radical $M( \lambda_0 - N\gamma ) \subset M( \lambda_0 )$ as 
we cross the reducibility hyperplane $H_{\gamma,N}$ (cf. \cite{Y}, 
Proposition 3.2).  Recall that $\varepsilon( H_{\gamma,N}, z )$ encodes 
information about singular vectors:  if $f \in U( \mfn^{op} )$ is such that 
$fv_{\lambda_0}$ generates $M( \lambda_0 -N\gamma ) = M( s_\gamma \lambda_0)$, 
then 
$$\sgn \left< fv_{\lambda_t-\rho}, fv_{\lambda_t-\rho} \right>_{\lambda_t} = 
\left \lbrace 
\begin{array}{ll}
\varepsilon( H_{\gamma,N}, z ) & \text{if }\lambda_t \in H_{\gamma,N}^+ \\
-\varepsilon( H_{\gamma,N}, z ) & \text{if } \lambda_t \in H_{\gamma,N}^-.
\end{array}
\right.
$$
Invariant Hermitian forms on Verma modules are unique up to a 
real scalar, which is determined by the inner product of a generator with 
itself.  The proposition now follows from the observation that $x\lambda = 
s_\alpha sx \lambda$ and $x\lambda = s_{x\alpha}xs \lambda$.
\end{proof}
\begin{corollary}
Letting $s=s_\alpha$ be a simple reflection, the signed Kazhdan-Lusztig 
polynomials satisfy:
\begin{itemize} 
\item[a)]
$P_{w_\lambda x, w_\lambda y}^{\lambda,w} = sgn( -w\rho, x\alpha)  
\varepsilon( H_{x\alpha,(xs\lambda,x\alpha^\vee)}, xs )
P_{w_\lambda xs, w_\lambda y}^{\lambda,w}$
if $ys > y$ and $x, xs \geq y$
\item[a')]
$P_{w_\lambda x, w_\lambda y}^{\lambda,w} = sgn( -w\rho, \alpha)  
\varepsilon( H_{\alpha,(sx\lambda,\alpha^\vee)}, sx )
P_{w_\lambda sx, w_\lambda y}^{\lambda,w}$
if $sy > y$ and $x, sx \geq y$.
\end{itemize}
\end{corollary}

An excellent companion for the remainder of this section is \cite{GJ} from 
which the results of this section are derived.  
The objective is to compute $P_{w_\lambda x, w_\lambda y}
^{\lambda,w}$ for case b).
We begin by introducing some background material.

\subsection{Category $\cO$} 
\label{category_O}
Bernstein-Gelfand-Gelfand defined category $\cO$ in \cite{BGG}.  It is the 
subcategory of the category of $\g$-modules consisting of modules $M$ 
satisfying:
\begin{enumerate}
\item $\displaystyle{M = \oplus_{\mu \in \h^*} M_\mu}$
\item $M$ is finitely generated
\item $M$ is $\mfn$-finite (i.e.\ $U(\mfn)v$ is finite-dimensional for 
every $v \in M$).
\end{enumerate}
Category $\cO$ is closed under arbitrary direct sums, quotients, submodules, 
and tensoring with finite-dimensional modules.  Verma modules are objects in 
category $\cO$ and the simple objects of $\cO$ consist of the irreducible 
highest weight modules $L( \mu )$ where $\mu \in \h^*$.   Irreducible highest
weight modules form an additive basis of the Grothendieck group of category 
$\cO$.

$\h^*$ is a disjoint union of $W$-orbits which are called {\bf blocks}.
Recall that $\chi_\mu = \chi_\nu$ if and only if $\nu \in W\mu$.  For each 
block $D$ and some $\mu \in D$, we define
$$\cO_D := \cO_\mu := \lbrace \text{modules}\in \cO \,|\, \exists \, N 
\text{ such that }
(z - \chi_\mu)^N \text{ annihilates } M \, \forall \, z \in Z( \g ) \rbrace.$$

\begin{remark}
These blocks are larger than the standard ones in ring theory:  two irreducible 
modules belong to the same block if they admit a non-trivial extension in the 
category.  In the case of category $\cO$, this amounts to two irreducible 
highest weight modules having the same infinitesimal character and all of 
their weights differing by sums of roots.
\end{remark}

Category $\cO$ decomposes into {\bf blocks} $\cO_D$ {\bf of category $\cO$}:  
\begin{theorem}(cf. \cite{BGG}, property 4) of Section 3)
$$\cO = \bigoplus_{\text{blocks } D} \cO_D.$$
\end{theorem}
Denote projection onto $\cO_D$ (or $\cO_\mu$) by $Pr_D$ (resp.\ $Pr_\mu$).  
Projection onto the blocks of
category $\cO$ defines what is known as the {\bf primary decomposition} 
$$M \simeq \bigoplus_{\text{blocks} D} Pr_D M$$
of a module $M$ in $\cO$.  For any $M \in \cO$, $Pr_D M$ is non-zero 
for finitely many $D$.  $Pr_D M$ is called the {\bf primary component of 
$M$ with respect to the block $D$}.
\begin{theorem}
\label{primary_decomp}
(cf. \cite{J2}, Satz 1, iv) of Section 3.)  Primary 
decomposition of a module in category $\cO$ which admits an invariant 
Hermitian form is an orthogonal decomposition into submodule pairs or 
singletons.  Specifically, 
$$Pr_D M \text{ and } Pr_{D'} M \text{ are orthogonal for } D' \neq 
- \bar{D}.$$
\end{theorem}
\begin{proof}
It is straightforward to modify Jantzen's proof that primary decomposition 
of modules in category $\cO$ admitting a contravariant form is an 
orthogonal decomposition.
\end{proof}

\subsection{Jantzen's translation functors and his determinant formula.}
\label{trans_functors}
For an integral weight $\mu$, let $F(\mu)$ denote the finite-dimensional 
representation of extremal weight $\mu$.
Jantzen's translation functors are 
compositions of tensoring with a finite-dimensional module with projections 
onto blocks of category $\cO$:
\begin{eqnarray*}
T_{\lambda}^{\lambda + \mu} : \lbrace \text{modules of inf'l 
character } \lambda \rbrace &\to& \lbrace \text{modules of inf'l 
character } \lambda + \mu \rbrace \\
M &\mapsto& Pr_{\lambda + \mu} ( M \otimes F(\mu) )
\end{eqnarray*}
for any $\lambda \in \h^*$ (cf. \cite{KO}).
We will later use $T_D^{D'}$ in place of $T_\lambda^{\lambda'}$ when $\lambda$ 
and $\lambda'$ are both antidominant, and therefore the value of $\lambda'
- \lambda$ may be recovered.

Recall for $\lambda \in \h^*$ that $\Delta_\lambda( \g, \h) = \lbrace \alpha
\in \Delta( \g, \h ) | \left< \lambda, \alpha^\vee \right> \in \bbZ \rbrace$.  
A {\bf facette} $\cF$ is a non-empty subset of some $\bbA( \lambda ) :=
\bbQ \Delta_\lambda( \g, \h) \otimes_\bbQ \bbR$ associated with a disjoint 
union $\Delta_\lambda^+( \g, \h) = \Delta_\cF^0 \cup \Delta_\cF^+ \cup 
\Delta_\cF^-$:
$$ \cF = \left \lbrace x \in \bbA( \lambda ) \left \vert
\begin{array}{rl}
\left< x, \alpha^\vee \right> = 0 & \text{if } \alpha \in \Delta_\lambda^0  \\
\left< x, \alpha^\vee \right> > 0 & \text{if } \alpha \in \Delta_\lambda^+  \\
\left< x, \alpha^\vee \right> < 0 & \text{if } \alpha \in \Delta_\lambda^-  
\end{array}
\right .
\right \rbrace.$$
If $M$ is irreducible, $\lambda$ and $\lambda + \mu$ are antidominant with 
$\lambda+\mu$ in the closure of the facette containing $\lambda$, then 
$T^{\lambda+\mu}_{\lambda}M$ is irreducible or zero (\cite{J}, Theorem 2.11).
If both $\lambda$ and $\lambda + \mu$ are strictly antidominant, then
$T^{\lambda + \mu}_{\lambda} : \cO_\lambda \to \cO_{\lambda+\mu}$ is 
an equivalence of categories.

We may extend the translation functor to category $\cO$:
\begin{eqnarray*}
T_{\lambda}^{\lambda + \mu} : \cO &\to& \cO \\
M &\mapsto& Pr_{\lambda + \mu} ( F(\mu) \otimes (Pr_\lambda M) ).
\end{eqnarray*}

In order to study how Jantzen's translation functors affect Verma modules and 
forms on Verma modules, we need some facts about the tensor product of a 
Verma module $M( \lambda )$ with a finite-dimensional module $F$:

\begin{theorem}(Bernstein-Gelfand-Gelfand, \cite{D}, Theorem 7.6.14.) 
Let $\mu_1$, $\ldots$, $\mu_N$ be an ordering of the weights of $F$ (with 
multiplicity) such that 
$\mu_i \leq \mu_j$ implies that $i \leq j$.  Then there is a filtration of 
$M = M( \lambda ) \otimes F$ by Verma modules:
$$ M = M^0 \supset M^1 \supset \cdots \supset M^N \supset M^{N+1} = \{ 0 \}$$
where $M^i / M^{i+1} \simeq M( \lambda + \mu_i )$.
\end{theorem}

\begin{theorem}(\cite{J2}, Satz 1, iii) of Section 3.)
\label{primary_generators} 
Let $M_i$ be $Pr_{\lambda + \mu_i}M$.  Then $M_i$ is generated as a 
$U( \mfn^{op} )$-module by the images of 
$v_{\lambda - \rho} \otimes F_\nu$ where $\nu$ is a weight of $F$ such that 
$\lambda + \nu$ belongs to the Weyl group orbit of $\lambda + \mu_i$.
\end{theorem}

Suppose the modules $U$ and $V$ admit an invariant Hermitian (resp. 
contravariant) form.  The tensor product of the two modules $U \otimes V$
naturally has an invariant Hermitian (resp. contravariant) form:
$\left< u_1 \otimes v_1, u_2 \otimes v_2 \right>_{U \otimes V} = 
\left< u_1 , u_2 \right>_U \cdot \left< v_1 , v_2 \right>_V.$
Since primary decomposition of a module is orthogonal with 
respect to invariant Hermitian forms in the case of a compact Cartan and also 
with respect to contravariant forms, the $M_i$s inherit invariant Hermitian 
forms from invariant Hermitian forms on $M( \lambda )$ and on $F$,  and they 
inherit contravariant forms from contravariant forms on $M( \lambda )$ and 
on $F$.  Jantzen has a determinant formula for such contravariant forms:
\begin{theorem} (\cite{J2}, Section 5.) 
\label{J_det_formula}
Suppose the numbers $(\lambda, \alpha^\vee )$ for $\alpha \in \Pi$ are 
algebraically independent over $\mathbb{Q}$.  Suppose $F = L( \lambda_0 )$ 
where $\lambda_0 \in \Lambda$ is strictly dominant and let $v_0$ be a 
highest weight vector of $F$
so that $(v_0, v_0) =1$.  Let $n(\mu)$ denote the multiplicity of the weight 
$\mu$ in $F$ and let $\lbrace e_{\mu,j} \rbrace_{1 \leq j \leq n(\mu)}$ be a 
$\bbZ$-basis of the $\mu$ weight space of $U( \mfn^{op} )_{\bbZ} v_0$.  $\mu =
\mu_i$ for some $i$.  Denote by $f_{\mu,j}$ the orthogonal projection of 
$v_{\lambda-\rho}\otimes e_{\mu,j}$ onto $M_i$.  The determinant of the 
contravariant form for the $f_{\mu,j}$ is $D_F( \mu )a_\mu$, where 
$D_F( \mu )$ is the determinant of the contravariant form with respect to a
$\bbZ$-basis of $U( \mfn^{op} )_{\bbZ}v_0$, and thus for the $e_{\mu,j}$, and 
$$a_\mu = \prod_{\alpha \in \Delta^+( \g,\h)} \,\, \prod_{
r > 0, r + \left< \mu, \alpha^\vee\right> \geq 0 }
\left( 
\frac{ (\lambda, \alpha^\vee ) - r}{(\lambda + \mu, \alpha^\vee ) + r} 
\right)^{n(\mu+r\alpha)}.$$
\end{theorem}

We now compare the canonical contravariant form with the Shapovalov form for 
the purpose of stating Jantzen's determinant formula for invariant 
Hermitian forms. 
We begin by introducing a $\bbZ_2$-grading of $\Lambda_r$ in the case of a 
compact Cartan.  From 
$\left[ \mfk, \mfk \right] \subset \mfk$,
$\left[ \mfk, \mfp \right] \subset \mfp$,
$\left[ \mfp, \mfp \right] \subset \mfk$, and from 
$\left[ \g_\alpha, \g_\beta \right] \subset \g_{\alpha+\beta}$, we see that 
a root is non-compact if and only if when expressed as a sum of simple roots, 
there are an odd number of non-compact roots in the sum, counting multiplicity.
It follows that for any $\mu \in \Lambda_r$, the parity of the number of 
non-compact roots in any expression of 
$\mu$ as a sum of roots is independent of the expression chosen.  We will 
denote the grading defined by this parity by $\varepsilon : \Lambda_r \to 
\bbZ_2$.

The Chevalley basis may be chosen so that 
$$X_\alpha^* = -\bar{X}_\alpha = \theta Y_\alpha =
(-1)^{\varepsilon( \alpha )} Y_\alpha = 
(-1)^{\varepsilon( \alpha )} \sigma( X_\alpha )$$
(cf. \cite{Y}).  Thus 
\begin{lemma} 
\label{form_comparison}
If $\h$ is compact and we choose a $\bbZ$-basis from 
$U( \mfn^{op} )_\bbZ v_{\lambda-\rho}$ for $(M( \lambda )_\bbZ)_\mu$, matrices 
representing the 
canonical contravariant form and the invariant Hermitian form with respect 
to this basis differ by multiplication by the scalar $(-1)^{\varepsilon(\lambda
-\rho -\mu)}$.
\end{lemma}
\begin{proposition}
\label{form_comp}
When $\h$ is compact,
Theorem \ref{J_det_formula} holds with ``invariant Hermitian form'' in place 
of ``contravariant form'' and 
$(-1)^{\varepsilon( \lambda_0 -\rho - \mu )n(\mu)}a_\mu$ 
in place of $a_\mu$.
\end{proposition}
\begin{proof}
A vector of weight $\lambda - \rho + \mu$ in $M( \lambda ) \otimes L( 
\lambda_0 )$ must be the sum of tensor products of a vector of weight 
$\lambda - \rho
-\nu$ and a vector of weight $\mu + \nu$ for some $\nu \in \Lambda_r$.  The 
proposition now follows from the grading, the lemma and the observation that 
$(-1)^{\varepsilon( \nu ) + \varepsilon( \lambda_0 - \rho - (\mu + \nu ) )}
= (-1)^{\varepsilon( \lambda_0 -\rho - \mu )}$.
\end{proof}

\subsection{Gabber and Joseph's generalization of category $\cO$}
For the purpose of studying the Kazhdan-Lusztig Conjecture,
Gabber and Joseph introduced modifications of category $\cO$.
Let $C \subset \h_A^*$ be of the form $\lambda + \Lambda_r$ (recall the 
discussion of the setup of \cite{GJ} after Theorem \ref{jantzen_filtration}).
Let $K_C$ be the subcategory of $U( \g_A )$-modules $M$ such 
that:
\begin{itemize}
\item[1)] $M = \sum_{\mu \in C - \rho} M_\mu$
\item[2)] $M$ is $U( \mfn_A )$-finite
\item[3)] $M$ is finitely generated over $U( \g_A )$.
\end{itemize}
Note that $M( \lambda )_A$ belongs to $\ob K_C$.

If $\mu \in C$, given any maximal ideal $m$ of $A$, there is a unique maximal 
submodule of $M( \mu )_A$ containing $mM( \mu )_A$ (cf.\ 1.7.2 of 
\cite{GJ}).  Call the corresponding simple quotient 
$L( m, \mu )$.  In our case, $A$ is a local ring, so we use $L( \mu )_A$ in 
place of $L( (t), \mu )$.

A {\bf block} in the context of the category $K_C$ is a subset $D$ of $C$ 
whose specialization at $t=0$, $\bar{D} = \{ \bar{
\lambda} \, | \, \lambda \in D \}$, is a $W$-orbit.  Define $$J_D = \bigcap_{
\mu \in D} \ker \chi_\mu.$$  
We may define, as we did for category $\cO$, the {\bf primary component of 
$M \in \ob K_C$ with respect to the block $D$}: 
$$Pr_D M = \{ m \in M \, | \, \text{for all } z \in J_D, \text{ exists } n \in 
\bbZ^+ \text{ such that } z^n m = 0 \}.$$

We note that $C = \coprod_i D_i$ is a countable union of blocks.  As for 
category $\cO$:

\begin{proposition} (cf. Proposition 1.8.4, \cite{GJ})
For $M \in \ob K_C$, we have 
$$M = \oplus_{i} Pr_{D_i} M,$$
the {\bf primary decomposition} of $M$.
\end{proposition}

Given a block $D \subset C$, the subcategory $K_D$ of $K_C$ consists of 
modules whose simple quotients are among the $L( m, \mu )$ where $m$ is a
maximal ideal of $A$ and $\mu \in D$.  $Pr_D$ takes objects in $K_C$ to 
objects in $K_D$.

In \cite{GJ}, Gabber and Joseph extended Jantzen's definition of translation 
functors to category $K_C$.  Let $D = W_\lambda \lambda + \delta t$, 
where $\lambda, \delta \in \h^*$ are regular and $\lambda$ is antidominant.
Let $D' = W_\lambda(\lambda - \mu) + \delta t$, where $\mu \in \Lambda$ and 
$\lambda - \mu$ is antidominant. 
$$ T_D^{D'} M = Pr_{D'} ( F( -\mu )_A \otimes_A ( Pr_D M ) )$$
is the translation functor from the block $K_D$ to the block $K_{D'}$.

We refer the reader to Definition 2.3 of \cite{Y} for the definition of the 
Hermitian dual of a module.  Given a module $M$ in $\ob K_{\bar{C}}$, we define 
$\delta^h(M)$ to be the $\h$-finite part of its Hermitian dual $M^h$.
\begin{lemma} 
\label{hermitian_dual}
Let $M \in \ob K_{\bar{D}}$ for some block $D \subset C$.
If $M$ admits a non-degenerate invariant Hermitian form and $\h$ is compact, 
then $\delta^h( M ) \cong M$.
\end{lemma}
\begin{proof}
We may modify Section 3.10 and Lemma 4.7 (iii) of \cite{GJ}.
\end{proof}

\subsection{Coherent continuation functors}
\label{CoherentContinuation_subsection}
Suppose $\lambda \in \h^*$ is antidominant and 
regular and $\delta \in \h^*$ is regular.  Let $D = W_\lambda \lambda + 
\delta t$.  Let $s = s_\alpha$ be a simple reflection in $W_\lambda$.  
We may choose $\nu_\alpha \in \Lambda$ so that $\lambda - 
\nu_\alpha$ is 
antidominant and so that the only root $\beta$ for which $( \lambda - 
\nu_\alpha, \beta ) = 0$ is $\beta = \alpha$.  Let $D_\alpha = W_\lambda( 
\lambda - \nu_\alpha) + \delta t$.  Since we are studying invariant 
Hermitian forms, we assume furthermore that $\delta$ and $\lambda$ are 
imaginary, although the statements which follow hold for non-imaginary 
$\delta$ and $\lambda$ if they contain no reference to invariant Hermitian 
forms.  {\it We fix this notation for the remainder 
of this article}.   The generalized notion (it exists for category $\cO$ 
also) of {\bf translation to the $\alpha$ wall} is the functor
$$T_D^{D_\alpha} M = Pr_{D_\alpha}
(F( -\nu_\alpha)_A \otimes_A ( Pr_D M ) )$$
and
$$T_{D_\alpha}^D M = Pr_D ( F( \nu_\alpha)_A \otimes_A (
Pr_{D_\alpha} M ) )$$
is {\bf translation from the $\alpha$ wall}.
Translation to the $\alpha$ wall followed by translation from the $\alpha$ 
wall, denoted by $\theta_\alpha = T_{D_\alpha}^D T_D^{D_\alpha}$,
is an exact functor known as {\bf coherent continuation across the 
$\alpha$ wall} or the {\bf reflection functor across the $\alpha$ wall}.
We will also use $\theta_\alpha$ to denote coherent continuation in category 
$\cO$.
Due to results in sections \ref{category_O} and \ref{trans_functors}, 
if $M$ carries an invariant Hermitian form, then so 
do $T_D^{D_\alpha} M$, $T_{D_\alpha}^D M$, and $\theta_\alpha M$ naturally.
We use $T \left< \cdot, \cdot \right>$ to denote the form which results from 
application of the translation functor $T$ to a module with invariant 
Hermitian form $\left< \cdot, \cdot \right>$.

We would like to describe the form which arises from a translation functor 
when the module and its invariant Hermitian form are a Verma module and its 
Shapovalov form respectively.  
For any $z \in W_\lambda$,
$$T_D^{D_\alpha}M( z\lambda + \delta t )_A \simeq M( z(\lambda-\nu_\alpha) + 
\delta t )_A$$
by Satz 2.9 of \cite{J}. 
We define
\begin{equation}
\label{gen_trans_to_wall}
v_{z(\lambda-\nu_\alpha) + \delta t - \rho}' := 
Pr_{D_\alpha} v_{z\lambda + \delta t - \rho} \otimes 
e_{-z\nu_\alpha,1},
\end{equation}
$$\left< \cdot , \cdot \right>_{z(\lambda - \nu_\alpha) + \delta t}' := 
T_D^{D_\alpha} \left< \cdot, \cdot \right>_{z\lambda + \delta t}, 
\qquad \text{and}$$
\begin{equation}
c_z' := \left< v_{z(\lambda-\nu_\alpha) + \delta t - \rho}',
v_{z(\lambda-\nu_\alpha) + \delta t - \rho}' \right>_{z(\lambda - \nu_\alpha)
+ \delta t}'.
\end{equation}
Let $\lambda_\alpha^- \in \Lambda^+$ be the highest weight of 
$F( -\nu_\alpha )$.
According to Theorems \ref{primary_generators} and \ref{J_det_formula} and 
Lemma \ref{form_comparison}, 
the form $\left< \cdot, \cdot \right>_{z(\lambda - \nu_\alpha) + \delta t}'$ 
on $M(  z(\lambda-\nu_\alpha) + \delta t)_A$ is such that
\begin{equation*}
c_z' =(-1)^{\varepsilon(\lambda_\alpha^- +
z\nu_\alpha)}D_{F(-\nu_\alpha)}(-z\nu_\alpha) 
a_{-z\nu_\alpha}'
\end{equation*}
where
\begin{equation}
\label{a'}
a_{-z\nu_\alpha}' = \prod_{\beta \in \Delta^+( \g,\h)} \,\, 
\prod_{ r > 0, r + \left< -z\nu_\alpha, \beta^\vee\right> \geq 0 }
\left( 
\frac{ (z\lambda + \delta t, \beta^\vee ) - r}{(z\lambda - z\nu_\alpha + 
\delta t, \beta^\vee ) + r} 
\right)^{n(-z\nu_\alpha+r\beta)}.
\end{equation}
On the level of signature characters, we have 
\begin{equation}
ch_s \overline{\left< \cdot, \cdot \right>}_{z(\lambda-\nu_\alpha)+\delta t}' = 
\sgn( c_z' ) ch_s \left< \cdot, \cdot \right>_{z(\lambda-\nu_\alpha)}.
\end{equation}
Here, we observe that Jantzen's determinant formula holds in the 
category $K_C$ setting also since we work with $U(\mfn^{op})_\bbZ$ bases and
hence his projection formulas and recursive formulas hold (cf. \cite{J2}, 
Section 5).

Returning to the problem of developing a
recursive formula for signed Kazhdan-Lusztig polynomials in case b), {\it we 
fix $x,y \in W_\lambda$ such that $y > ys$ and $x < xs$ for the remainder of 
this section}.
Following the notation of 
\cite{GJ}, let $X = M( xs \lambda + \delta t )_A$ and let $Z = M( x \lambda 
+ \delta t)_A$.  Define $Y$ to be $\theta_\alpha Z$.  Then:
\begin{proposition} (\cite{GJ}, section 3.6.)
\label{Ystructure}
\begin{itemize}
\item[i)] $\theta_\alpha X \simeq \theta_\alpha Z \simeq Y$.  
\item[ii)] 
There is a short exact sequence $0 \to X \to Y \stackrel{\pi}{\to} 
Z \to 0$.
\end{itemize}
\end{proposition}
\begin{remark}
\label{extension}
Because of the short exact sequence, $Y$ is called an {\bf extension of $X$ 
by $Z$}.
\end{remark}
\begin{remark}
Gabber and Joseph's results are for contravariant forms.  The invariant 
Hermitian form analogues of their results hold in the compact Cartan case:
their proofs may be transferred to the Hermitian form setting 
using Lemma \ref{Jantzen_same} and Theorem \ref{primary_decomp}.
\end{remark}

We review Gabber and Joseph's discussion of the filtration of $Y$ by Verma 
modules.  Now for 
$x' \in W_\lambda$, 
\begin{equation}
\label{pfiltration}
\left[ T_{D_\alpha}^D M( x( \lambda - \nu_\alpha ) + \delta t )_A 
: M( x'\lambda + \delta t )_A \right] = \dim (F(\nu_\alpha)_A)_{\mu}
\end{equation}
where $\mu = x'\lambda - x(\lambda - \nu_\alpha) = x'\lambda - xs( \lambda
- \nu_\alpha )$.  According to Satz 2.9 of 
\cite{J}, there are two solutions:  $x' = x, \mu = x\nu_\alpha$ and 
$x' = xs, \mu = xs \nu_\alpha$.  By Theorem 
\ref{primary_generators}, $Y$ is generated as a $U( \mfn^{op} )$ module by 
$v_{x\lambda+\delta t-\rho}'' = Pr_{D_\alpha} \left( v_{x(\lambda-\nu_\alpha) 
+ \delta t -\rho}' \otimes e_{x\nu_\alpha,1} \right)$ 
and by $v_{xs\lambda+\delta t -\rho}'' =
Pr_{D_\alpha} \left( v_{x(\lambda-\nu_\alpha) + \delta t -\rho}' \otimes e_{x
s\nu_\alpha,1} \right)$.
Observe that $v_{x\lambda+\delta t-\rho}''$ and  $v_{xs\lambda+\delta t
-\rho}''$ are mutually orthogonal with respect to
$\left< \cdot, \cdot \right>_D'' := T^D_{D_\alpha} \left< \cdot, \cdot 
\right>_{x(\lambda-\nu_\alpha) + \delta t}'$.

Recall $c_x' = \left< v_{x(\lambda-\nu_\alpha) + \delta t -\rho}' ,
v_{x(\lambda-\nu_\alpha) + \delta t -\rho}' \right>_{x(\lambda-\nu_\alpha) + 
\delta t}'$.
Let $\lambda_\alpha^+ \in \Lambda$ be the highest weight of 
$L( \lambda_\alpha^+ )_A = F( \nu_\alpha )_A$.
From Theorem \ref{J_det_formula}, we have
\begin{eqnarray*}
\left< v_{x\lambda + \delta t - \rho}'' , 
v_{x\lambda + \delta t - \rho}'' \right>_D''
 &=&(-1)^{\varepsilon(\lambda_\alpha^+ - x\nu_\alpha)}
D_{F(\nu_\alpha)_A}( x\nu_\alpha) a_{x\nu_\alpha}'' c_x'
\qquad \text{and} \\
\left< v_{xs\lambda + \delta t - \rho}'' , 
v_{xs\lambda + \delta t - \rho}'' \right>_D''
&=& (-1)^{\varepsilon(\lambda_\alpha^+ - xs\nu_\alpha)} 
D_{F(\nu_\alpha)_A}( xs\nu_\alpha) a_{xs\nu_\alpha}''  c_x'
\end{eqnarray*}
where
\begin{eqnarray*}
a_{x\nu_\alpha}'' &=&  \prod_{\beta \in \Delta^+} \,\prod_{r > 0, r + \left<
x\nu_\alpha, \beta^\vee \right> \geq 0} \left( \frac{( x(\lambda-\nu_\alpha)
+ \delta t, \beta^\vee) -r}{(x\lambda + \delta t, \beta^\vee) + r} \right)
^{n( x\nu_\alpha + r\beta)} 
\qquad\text{and} \\
a_{xs\nu_\alpha}'' &=& \prod_{\beta \in \Delta^+} \,\prod_{r > 0, r + \left<
xs\nu_\alpha, \beta^\vee \right> \geq 0} \left( \frac{( xs(\lambda-\nu_\alpha)
+ \delta t, \beta^\vee) -r}{(xs\lambda + \delta t, \beta^\vee) + r} \right)
^{n( xs\nu_\alpha + r\beta)}.
\end{eqnarray*}
We compute which factors are zero at $t=0$.

\begin{list}{}{
\setlength{\itemindent}{0.3in}
\setlength{\leftmargin}{0.2in}}
\item[Denominator of $a_{x\nu_\alpha}''$:]
We require $r$ and $\beta > 0$ such that $( x\lambda, \beta^\vee ) = -r < 0$.  
Then $s_\beta x \lambda - x \lambda = r\beta$ so $s_\beta x \lambda -
x(\lambda - \nu_\alpha) = x\nu_\alpha + r\beta $.  From (\ref{pfiltration})
and from part ii) of Theorem \ref{Ystructure},
we see that 
$$n( x\nu_\alpha + r\beta ) = 
\left \lbrace \begin{array}{cl}
1 &\text{if } s_\beta x = xs_\alpha \Rightarrow \beta = x \alpha \\
0 &\text{otherwise.}
\end{array}
\right.
$$
We conclude that the denominator has exactly one factor, $(\delta t, 
x\alpha^\vee )$, which is zero at $t=0$.
\item[Numerator of $a_{x\nu_\alpha}''$:]
Suppose we have $\beta > 0$ and $( x( \lambda - \nu_\alpha ), \beta^\vee ) =
r > 0$.  Then $x( \lambda - \nu_\alpha ) - s_\beta x( \lambda - \nu_\alpha )
 = r\beta$ so $s_\beta x \nu_\alpha -r\beta = s_\beta x\lambda - x(\lambda - 
\nu_\alpha)$.  By (\ref{pfiltration}) and by part ii) of Theorem 
\ref{Ystructure},
$$n( x\nu_\alpha + r\beta ) = n( s_\beta x \nu_\alpha -r\beta ) = 
\left \lbrace \begin{array}{cl}
1 &\text{if } s_\beta x = xs_\alpha \Rightarrow \beta = x \alpha \\
0 &\text{otherwise.}
\end{array}
\right.
$$
However, $(x( \lambda - \nu_\alpha ) , x\alpha^\vee ) = 0 \neq r$ and we 
deduce that the numerator has no factors which are zero at $t=0$.
\end{list}
Similarly, none of the factors in the numerator and the denominator of 
$a_{xs\nu_\alpha}''$ are zero at $t=0$.

\begin{remark}
The results of this section hold with any $z < zs$ in place of $x$.  We 
define $v_{z\lambda + \delta t - \rho}''$, $v_{zs\lambda + \delta t - \rho}''$,
$a_{z\nu_\alpha}''$, 
and $a_{zs\nu_\alpha}''$ analogously for all such $z \in W_\lambda$.
\end{remark}

\subsection{A recursive formula in the difficult case.}
Here, we combine the results of the preceding subsections to deduce a 
recursive formula for computing signed Kazhdan-Lusztig polynomials for case 
b).

We will need Gabber and Joseph's description of $\bar{Y}_j$ where the 
form on $Y$ arises from the form on $Z$ and coherent continuation (cf. 
4.4, 4.5, 4.6 of \cite{GJ}).
Recall the exact 
sequence from Theorem \ref{Ystructure}.  If we define $\bar{Y}^x_j = 
( \bar{Y}^j \cap \bar{X} ) / (\bar{Y}^{j+1} \cap \bar{X} )$ and 
$\bar{Y}^z_j = \pi( \bar{Y}^j ) / \pi( \bar{Y}^{j+1} )$  then there is 
a short exact sequence
$$0 \to \bar{Y}^x_j \to \bar{Y}_j \to  \bar{Y}^z_j \to 0.$$
For $M \in \ob K_C$, we define $M^+$ (resp.\ $M^-$) to be the smallest 
(resp.\ largest) submodule of $M$ for which $\theta_\alpha(M / M^+ ) = 0$ 
(resp.\ $\theta_\alpha M^- = 0$ ).  We have the short exact sequences
$$ 0 \to \bar{X}^+_{j+1} \to \bar{Y}^x_j \to \bar{X}^-_j \to 0$$
and
$$ 0 \to \bar{Z}^-_{j+1} \to \bar{Y}^z_j \to \bar{Z}^+_j \to 0$$
(cf.\ 4.5 (2), 4.5 (4), Lemma 4.6 ii) and Proposition 4.7 of \cite{GJ}).  This 
gives us the four-step filtration of $\bar{Y}_j$:
\setlength{\extrarowheight}{2pt}
\begin{equation*}
\begin{array}{|c|c|c|c|}
\hline
\multicolumn{4}{|c|}{ \bar{Y}_j } \\ \hline
\multicolumn{2}{|c|}{ \bar{Y}^x_j} & 
\multicolumn{2}{|c|}{ \bar{Y}^z_j} \\ \hline
\bar{X}^+_{j+1} & \bar{X}^-_j & \bar{Z}^-_{j+1} & \bar{Z}^+_j  \\ \hline
\end{array}
\end{equation*}
Here we remark that in the paper \cite{GJ}, because 3.14 still holds and 
because we 
may modify Lemma 3.15 for invariant Hermitian forms, Hermitian analogues of 
results in sections 4.4 to 4.7 hold.  We have:
\begin{lemma} (cf. \cite{GJ}, Lemma 4.5.)

\begin{itemize}
\item[i)] $\left< \bar{X}^+_{j+1}, \bar{Y}^x_j \right> = 0$
\item[ii)] $\left< \bar{X}^+_{j+1}, \ker( \bar{Y}_j \to \bar{Z}^+_{j} )
\right> = 0$.
\end{itemize}
\end{lemma}
Furthermore, $\delta^h( \bar{X}^+_{j+1} ) \cong \bar{Z}^+_j$. 
It follows that:
\begin{proposition}
\label{Ysig}
Consider the four-step filtration of $\bar{Y}_j$.
$\bar{X}^+_{j+1}$ is paired with $\bar{Z}^+_j$, whence the signature 
character of $\bar{Y}_j$ is given by the signature characters of 
$\bar{X}^-_j$ and $\bar{Z}^-_{j+1}$.
\end{proposition}

We will clarify the latter half of this statement, which is vague. 
First, we discuss the structure of $\bar{X}^-_j$ and $\bar{Z}^-_{j+1}$, 
which is given (along with the structure of $\bar{X}_i^\pm$ and $\bar{Z}_i^\pm$ 
for any $i$) by Kazhdan-Lusztig polynomials and the following proposition:

\begin{proposition} (\cite{GJ}, Lemma 3.6, 3.11.)
\begin{itemize}
\item[i)] $\theta_\alpha L( z\lambda ) = 0$ if $z > zs$ and 
$\theta_\alpha L( z \lambda ) \neq 0$ otherwise.
\item[ii)] When $z < zs$, $\theta_\alpha L( z\lambda )$ has a unique simple 
quotient and it is isomorphic to $L( z\lambda )$.  
The corresponding unique maximal submodule has unique simple submodule
$L(z\lambda )$. 
\end{itemize}
\end{proposition}
Recalling that $\bar{X}_i$ is semisimple, $\bar{X}_i^+$ is a sum of simple 
submodules $L(z \lambda )$ for which $z > zs$ and $\bar{X}_i^-$ is a sum of
simple submodules $L( z \lambda )$ for which $z < zs$.  Likewise for 
$\bar{Z}_i$.

We study $\theta_\alpha L( z \lambda )$ when $z < zs$ in more detail.
The functor $T^{\bar{D}}_{\bar{D}_\alpha}$ is a left and a right adjoint to 
$T^{\bar{D}_\alpha}_{\bar{D}}$
(cf. 3.4 of \cite{GJ} and (3.5) of \cite{V2}).  It follows that
\begin{eqnarray*}
\Hom_\g\left( L( z\lambda ), T^{\bar{D}}_{\bar{D}_\alpha}T^{\bar{D}_\alpha}_{
\bar{D}} L( z\lambda ) \right)
&\cong& \Hom_\g\left( T^{\bar{D}_\alpha}_{\bar{D}} L( z\lambda ), 
T^{\bar{D}_\alpha }_{\bar{D}} L( z\lambda ) \right) \\
&\cong& \Hom_\g\left( T^{\bar{D}}_{\bar{D}_\alpha} T^{\bar{D}_\alpha}_{\bar{D}} 
L( z\lambda ),  L( z\lambda ) \right)
\end{eqnarray*}
from which we obtain a chain complex
$$0  \to L(z \lambda ) \stackrel{i}{\to} \theta_\alpha L(z\lambda ) 
\stackrel{p}{\to} L(z\lambda )
\to 0$$
(cf. \cite{V2}, Theorem 3.7).
Because $L( z\lambda )$ is simple, the first map is injective and the 
second map is surjective.

Since $\theta_\alpha L(z\lambda )$ admits a non-degenerate invariant 
Hermitian form (for example, the form acquired through coherent continuation
and the form on $L( z\lambda )$), 
$\delta^h ( \theta_\alpha L( z \lambda ) ) \cong
\theta_\alpha L( z \lambda )$.  $\delta^h$ takes submodules of a module $M$ 
to quotients of $\delta^h M$ and quotients of $M$ to submodules of 
$\delta^h M$.  Since $\delta^h$ does not take the submodule $L(z\lambda )$ of
$\theta_\alpha L(z\lambda)$ to the submodule $L(z\lambda )$ of $\delta^h
(\theta_\alpha L(z\lambda) )$, it follows that that submodule 
cannot be paired with itself, and hence it 
is paired with the quotient $L( z\lambda )$.  We conclude:

\begin{lemma}
\label{sig_from_Ualpha}
Suppose $z < zs$.  Then $ch_s \theta_\alpha L(z \lambda ) = ch_s U_\alpha
L( z\lambda )$ where $U_\alpha L(z \lambda )$ is defined to be the 
cohomology of the complex
$$ 0 \to L( z \lambda ) \hookrightarrow \theta_\alpha L( z \lambda ) 
\twoheadrightarrow L( z \lambda ) \to 0 .$$
\end{lemma}

$U_\alpha$ may be extended to semisimple modules via $U_\alpha (M \oplus N )
= U_\alpha M \oplus U_\alpha N$.  In particular, we may apply $U_\alpha$ to 
$\bar{Z}_j^+$.  Since $\theta_\alpha M^j = (\theta_\alpha M )^j$ for 
$M \in \ob K_C$ (cf.  \cite{GJ}, Lemma 4.3 ii) ), we see that $\bar{Y}_j
 = \theta_\alpha \bar{Z}_j = \theta_\alpha \bar{Z}_j^+$.
\begin{proposition} (cf. \cite{GJ}, Proposition 4.7 iv).) 
\label{Ualpha_SES}
There is a short exact sequence
$$ 0 \to \bar{X}_j^- \to U_\alpha \bar{Z}_j^+ \to \bar{Z}_{j+1}^- \to 0.$$
Furthermore, choosing the form $\theta_\alpha \left< \cdot, 
\cdot \right>_{x\lambda + \delta t}$ on $Y$:
$$ch_s \bar{Y}_j = ch_s U_\alpha \bar{Z}_j^+
= \sgn( \bar{c}_{xs}'' \bar{c}_x' )ch_s \bar{X}_j^- + 
\sgn( \bar{c}_x'' (\delta, x\alpha^\vee ) \bar{c}_x' )ch_s \bar{Z}_j^-  
$$
where for $z < zs$
\begin{align*}
&c_z'' := (-1)^{\varepsilon( \lambda_\alpha^+ 
-z\nu_\alpha)} D_{F(\nu_\alpha)}( z\nu_\alpha )a_{z\nu_\alpha}''
( \delta t, z \alpha^\vee) \\
\text{and} \qquad  &
c_{zs}'' := (-1)^{\varepsilon( \lambda_\alpha^+ -zs\nu_\alpha)}
D_{F(\nu_\alpha)}( zs\nu_\alpha )a_{zs\nu_\alpha}''.
\end{align*}
\end{proposition}
\begin{proof}
This follows from our previous discussion, Theorems \ref{primary_generators}
and \ref{J_det_formula}, and our analysis of $a_{z\nu_\alpha}''$ and 
$a_{zs\nu_\alpha}''$.
\end{proof}

We discuss the signature character of an invariant Hermitian 
form on some $U_\alpha L( z\lambda)$.  The process of coherent continuation 
in category $\cO$ uniquely determines an invariant Hermitian form on 
$U_\alpha L(z \lambda )$ from a form on $L(z \lambda )$.  Since $U_\alpha
L( z\lambda )$ is semisimple by Vogan's Conjecture, it may have many other 
non-degenerate invariant Hermitian forms.  
For example, another natural form on $U_\alpha 
L(x\lambda )$ is the form on $\bar{Y}_0$ 
which arises from the Jantzen filtration of 
$\theta_\alpha Z = \theta_\alpha M( x\lambda + \delta t )$. 
The signature depends on $\delta$ while the form coming from coherent 
continuation in category $\cO$ does not, and so the signatures may be 
different.  The form given by the Jantzen filtration is the form in which 
we are interested.  In the following computations, {\it we will always study 
the form arising from the Jantzen filtration.}

We compute the signature character of the form on 
$\overline{\left( \theta_\alpha M( z\lambda + \delta t) \right)}_0 \cong 
\theta_\alpha L( z \lambda )$ for $z < zs$.  By Proposition 
\ref{Ualpha_SES}, there is a short exact sequence
$$ 0 \to \overline{ M( zs \lambda + \delta t) }^-_0 \to
U_\alpha \overline{ M( z \lambda + \delta t) }_0 \to
\overline{ M( z \lambda + \delta t) }^-_1 \to 0.$$
By $ch_s U_\alpha \overline{ M( z \lambda + \delta t) }_0$ we mean the 
signature character of the form given by the Jantzen filtration in the 
direction $\delta$.
By our analysis
of the denominators and numerators of $a_{z\nu_\alpha}''$ and of
$a_{zs\nu_\alpha}''$, the denominators of $c_z''$ and $c_{zs}''$ do not vanish 
at $t=0$.
From Jantzen's determinant formula and our short exact sequence 
above,
\begin{align*}
ch_s U_\alpha \overline{ M( z \lambda + \delta t) }_0
= & \sgn( \bar{c}_{zs}''\bar{c}_z' )
ch_s L( zs \lambda )  \\
& + \sgn( \bar{c}_z'' ( \delta, z \alpha^\vee) \bar{c}_z' )
\sum_{y \in W_\lambda | y > ys} a^{z\lambda,w}_{y,1}
ch_s L( y\lambda ).
\end{align*}
Using this in conjunction with the previous proposition gives:
\begin{proposition} 
If $x,y \in W_\lambda$ are such that $x < xs$ and $y > ys$ and $x > y$ 
then:
\begin{align*}
\sgn( \bar{c}_{xs}'' \bar{c}_x' ) P^{\lambda,w}_{w_\lambda xs, w_\lambda y}(q)
&+ \sgn( \bar{c}_x'' ( \delta, x\alpha^\vee ) \bar{c}_x' )q P^{\lambda, w}
_{w_\lambda x, w_\lambda y}(q) \\
= &\sum_{z \in W_\lambda | z < zs} \sgn( \bar{c}_z'' (\delta, z\alpha^\vee) 
\bar{c}_z')
a^{z\lambda,w}_{y,1} q^{\frac{\ell(z)-\ell(y) + 1}{2}} P^{\lambda,w}
_{w_\lambda x, w_\lambda z} (q)  \\
&+ \sgn( \bar{c}_{ys}''( \delta, ys\alpha^\vee )\bar{c}_{ys}') P_{w_\lambda x, 
w_\lambda ys}^{\lambda,w}(q) .
\end{align*}
\end{proposition}

We discuss the values of $\sgn( \bar{c}_z'')$, $\sgn( \bar{c}_{zs} )$ and 
$\sgn(  \bar{c}_z')$ for $z < zs$.

\begin{lemma}
For an integral weight $\nu$ and for all $w \in W$,
$$\sgn \left( D_{F(\nu)}(w\nu ) \right)= 1$$
(see Theorem \ref{J_det_formula} for notation).
\end{lemma}
\begin{proof}
We prove this by induction on $\ell( w )$.  Clearly this is true for $w = 1$.
We may assume $\nu$ to be dominant and let $v_\nu$ be the canonical generator 
of $F(\nu)$.  Suppose the lemma holds for $w \in W$ and $s_\alpha$ is a simple 
reflection such 
that $s_\alpha w > w$.  Let $a \in U(\mfn^{op})$ be such that $av_\nu$ is a 
vector of weight $w\nu$ in $F(\nu)$.  Now $s_\alpha w > w$, so $(w\nu, 
\alpha^\vee ) > 0$.  Let $n_\alpha = ( w\nu, \alpha^\vee ) = (\nu, 
w^{-1}\alpha^\vee ) \in \bbZ^{\geq 0}$.  $s_\alpha w\nu = w\nu - n_\alpha 
\alpha$ and so $Y_\alpha^{n_\alpha}av_\nu$ is a vector of weight $s_\alpha
w\nu$ in $F(\nu )$.  Because $av_\nu$ is a vector of extremal weight $w\nu$ 
and because the set of weights of $F(\nu)$ is convex, $Y_\alpha^{n_\alpha}
av_\nu \in F(\nu)$ implies $X_\alpha av_\nu = 0$.  Therefore
\begin{eqnarray*}
( Y_\alpha^{n_\alpha} av_\nu, Y_\alpha^{n_\alpha} av_\nu ) &=& 
( \sigma( Y_\alpha^{n_\alpha} a )Y_\alpha^{n_\alpha} a v_\nu, v_\nu)  
= ( \sigma( a ) X_\alpha^{n_\alpha} Y_\alpha^{n_\alpha} a v_\nu, v_\nu ) \\
(\text{from } X_\alpha a v_\nu = 0 )
&=& ( \sigma( a ) p( X_\alpha^{n_\alpha}Y_\alpha^{n_\alpha} ) a v_\nu, 
v_\nu ) \\
(\text{from }\mathfrak{sl}_2 \text{ theory})
&=& w\nu \left( H_\alpha \left( H_\alpha - 1  \right) \cdots \left( H_\alpha - 
(n_\alpha - 1 ) \right) \right) \cdot ( \sigma( a ) a v_\nu, v_\nu ).
\end{eqnarray*}
Now $w\nu \left( H_\alpha \left( H_\alpha - 1  \right) \cdots \left( H_\alpha - 
(n_\alpha - 1 ) \right) \right) > 0$  since $w\nu( H_\alpha ) = n_\alpha$.
By our induction hypothesis, $( \sigma( a ) a v_\nu, v_\nu ) > 0$.
Thus $( Y_\alpha^{n_\alpha} av_\nu, Y_\alpha^{n_\alpha} av_\nu ) > 0$, 
proving our lemma.
\end{proof}

\begin{remark}
We may also prove the lemma using the following unpublished result of Birgit 
Speh:  if $\h$ is a compact Cartan subalgebra, then given the finite 
dimensional representation of highest weight $\lambda_0$, the Shapovalov 
form is definite on each weight space, with the form being positive definite
(resp.\ negative definite) on the $\lambda_0-\mu$ weight space if 
$\varepsilon( \mu ) = 0$ (resp. $\varepsilon( \mu) = 1$).  Again, we may 
take $\nu$ to be dominant.  Comparing what Speh's formula and Lemma 
\ref{form_comparison} imply for the signature of the one-dimensional weight 
space corresponding to $w\nu$, we have $(-1)^{\varepsilon( \nu-w\nu)} = 
\sgn( D_{F(\nu)}(w\nu))(-1)^{\varepsilon( \nu - w\nu)}$ from which the 
lemma follows.
\end{remark}

\begin{lemma}  For $z < zs \in W_\lambda$:
\begin{itemize}
\item[i)] $\sgn( \bar{a}_{-z\nu_\alpha}' ) = 1.$
\item[ii)] $\sgn( \overline{ a_{z\nu_\alpha}''( \delta t, z\alpha^\vee)} ) = 
-1$.
\item[iii)] $\sgn{\bar{a}_{zs\nu_\alpha}'' }= 1$.
\end{itemize}
\end{lemma}
\begin{proof}
i):  Consider (\ref{a'}).  If $( -z\nu_\alpha, \beta^\vee ) \geq 0$ then 
$n( -z\nu_\alpha + r\beta ) = 0$ for $r > 0$.
If $( -z\nu_\alpha, \beta^\vee ) < 0$, then the index for the second product 
starts at $r = -( -z\nu_\alpha, \beta^\vee )$.
$$-z\nu_\alpha - (-z\nu_\alpha, \beta^\vee) \beta = s_\beta( -z\nu_\alpha )$$
is an extremal weight of $F( -\nu_\alpha )$.  Therefore $n( -z\nu_\alpha + 
r\beta ) = 0$ for $r > -(-z\nu_\alpha, \beta^\vee )$.  Therefore 
$\bar{a}_ {-z\nu_\alpha}'$ may be written
\begin{eqnarray*}
\bar{a}_{-z\nu_\alpha}' &=& 
\prod_{\beta \in \Delta^+( \g,\h), (-z\nu_\alpha, \beta^\vee ) < 0}
\left( \frac{(z\lambda , \beta^\vee) + (-z\nu_\alpha, \beta^\vee) }
{(z\lambda -z\nu_\alpha, \beta^\vee) - (-z\nu_\alpha,\beta^\vee)}
\right)^1 \\
&=& \prod_{\beta \in \Delta^+( \g,\h), (-z\nu_\alpha, \beta^\vee ) < 0}
\frac{(z(\lambda-\nu_\alpha) , \beta^\vee) }
{(z\lambda , \beta^\vee)}.
\end{eqnarray*}
$\lambda-\nu_\alpha$ lies in the closure of the antidominant Weyl chamber, 
which is the Weyl chamber to which $\lambda$ belongs.  Since $(\lambda-
\nu_\alpha, \beta^\vee) \neq 0$ for $\beta \neq \alpha$, we conclude that 
$\sgn( z( \lambda - \nu_\alpha), \beta^\vee ) = \sgn( z\lambda, \beta^\vee )$ 
for $\beta \neq z\alpha$.  Observing that $( -z\nu_\alpha, z\alpha^\vee ) > 0$,
we conclude that $\sgn( \bar{a}_{-z\nu_\alpha}') = 1$.

ii):  As in the previous case,
$$\overline{ a_{z\nu_\alpha}''( \delta t, z\alpha^\vee)} = 
( z \lambda, z \alpha^\vee )
\prod_{\beta \in \Delta^+( \g, \h) \setminus \{ z\alpha \}, 
( z\nu_\alpha, \beta^\vee ) < 0}
\left( \frac{ ( z\lambda, \beta^\vee ) }
{(z( \lambda - \nu_\alpha), \beta^\vee )} \right)^1$$
so $\sgn( \overline{a_{z\nu_\alpha}''( \delta t, z \alpha^\vee)} ) = -1$.

iii):  As in the first case, 
$$ \bar{a}_{zs\nu_\alpha}'' =
\prod_{\beta \in \Delta^+( \g, \h), 
( zs\nu_\alpha, \beta^\vee ) < 0}
\frac{( zs\lambda, \beta^\vee)}
{(zs( \lambda - \nu_\alpha), \beta^\vee ) }.
$$
Since $( zs\nu_\alpha, z\alpha^\vee ) > 0$, we conclude that $\sgn(
\bar{a}_{zs\nu_\alpha}'') = 1$.
\end{proof}

Combining the results of this subsection, cancelling out common factors of 
$(-1)^{\varepsilon( \lambda_\alpha^+ + \lambda_\alpha^-)}$, and observing 
that $x\nu_\alpha - xs\nu_\alpha = x( (\nu_\alpha, \alpha^\vee)\alpha )
= x( ( \lambda, \alpha^\vee )\alpha )$,
we arrive at:

\begin{theorem}
\label{signedKL_formula}
Letting $s = s_\alpha$ be a simple reflection, the signed Kazhdan-Lusztig 
polynomials are defined by the intial conditions $P_{x,x}^{\lambda,w} = 1$, 
$P_{x,y}^{\lambda,w} = 0$ for $x > y$ and the recursive formulas:
\begin{itemize}
\item[a)]
$P_{w_\lambda x, w_\lambda y}^{\lambda,w} = sgn( -w\rho, x\alpha)  
\varepsilon( H_{x\alpha,-(\lambda,\alpha^\vee)}, xs )
P_{w_\lambda xs, w_\lambda y}^{\lambda,w}$
if $ys > y$ and $xs > x \geq y$
\item[a')]
$P_{w_\lambda x, w_\lambda y}^{\lambda,w} = sgn( -w\rho, \alpha)  
\varepsilon( H_{\alpha,(sx\lambda,\alpha^\vee)}, sx )
P_{w_\lambda sx, w_\lambda y}^{\lambda,w}$
if $sy > y$ and $sx > x \geq y$
\item[b)]
If $x,y \in W_\lambda$ are such that $x < xs$ and $y > ys$ and $x > y$ 
then:
\begin{align*}
&-(-1)^{\varepsilon( ( \lambda, \alpha^\vee )x\alpha )} P^{\lambda,w}_{w_\lambda xs, 
w_\lambda y}(q)
+ \sgn( \delta, x\alpha^\vee ) q P^{\lambda, w}
_{w_\lambda x, w_\lambda y}(q) \\
&= \sum_{z \in W_\lambda | z < zs} \sgn( \delta, z\alpha^\vee) 
a^{z\lambda,w}_{y,1} q^{\frac{\ell(z)-\ell(y) + 1}{2}} P^{\lambda,w}
_{w_\lambda x, w_\lambda z} (q)  
+ \sgn( \delta, ys\alpha^\vee ) P_{w_\lambda x, 
w_\lambda ys}^{\lambda,w}(q) .
\end{align*}
\end{itemize}
\end{theorem}

\section{Some Examples}
\label{Examples_section}
\subsection*{Example 1: $\g_0 = \mathfrak{so}(2)$.}
We have $\h = \mft$.  Let $\Delta^+( \g, \h ) = \lbrace \alpha_1 \rbrace$ 
and let $\lambda_1$ be the corresponding fundamental weight.

Irreducible Verma modules:  Choose $\lambda \in \h^*$ so that $(\lambda, 
\alpha_1^\vee) \in (n,n+1)$ where $n \in \bbZ_{\geq 0}$.  Then $\lambda \in
A( n\lambda_1, w_0)$.  The reducibility hyperplanes separating the alcove 
$aA_0$ containing $\lambda$ and $\widetilde{a}A_0$ are $H_{\alpha_1, 1}$, 
$H_{\alpha_1, 2}$, $\ldots$ $H_{\alpha_1, n}$.  In the setup of Theorem 
\ref{Verma_formula} we choose the path so that $r_1 = s_{\alpha_1,n}$, 
$r_2 = s_{\alpha_1, n-1}$, $\ldots$, $r_n = s_{\alpha_1, 1}$.  Suppose 
$S \subset \lbrace 1, 2, \ldots, n \rbrace$ and $|S| \geq 2$.  Then 
$\overline{r_{i_1}}C_{i_2-1}$ and $\overline{r_{i_1}}C_{i_2}$ lie in the 
Wallach region, and thus  $\varepsilon( \overline{r_{i_1}}C_{i_2-1}, 
\overline{r_{i_1}}C_{i_2} ) = 0$.  Therefore $\varepsilon( S ) = 0$ for 
$|S| \geq 2$.  For our choice of path, note that
$C_i \supset (n-i, n-i+1)$, whence
$\varepsilon( \lbrace i \rbrace ) = \varepsilon( C_{i-1}, C_i ) = 
\varepsilon( H_{\alpha_1, n-i+1}, s_1) = \delta_{\alpha_1}^{n-i+1} =1$
(see Definition 5.2.16 and Lemma 5.2.17 or Theorem 6.12 of \cite{Y}).  
Substituting these values into 
Theorem \ref{Verma_formula}:
\begin{eqnarray*}
R^{A(n\lambda_1, w_0)} = ch_s M( \lambda ) &=& 
\frac{\sum_{i=1}^n 2 e^{\overline{r_i}r_i\lambda-\rho} + e^{\lambda-\rho}}
{\displaystyle{
\prod_{\alpha \in \Delta^+( \mfp, \mft)}( 1 - e^{-\alpha})
\prod_{\alpha \in \Delta^+( \mfk, \mft)}( 1 - e^{-\alpha})}} \\
&=& 
\frac{\sum_{i=1}^n 2 e^{\lambda - i \alpha_1 -\rho} + e^{\lambda-\rho}}
{1+e^{-\alpha_1}} \\
&=& 
\frac{\sum_{i=1}^n e^{\lambda - (i-1) \alpha_1 -\rho} + e^{\lambda - i 
\alpha_1 -\rho} } {1+e^{-\alpha_1}} \\
&=& e^{\lambda-\rho} + e^{\lambda-\rho-\alpha_1 -\rho} + \cdots +
e^{\lambda-(n-1)\alpha_1 -\rho} + \frac{e^{\lambda-n\alpha_1-\rho}}
{1+e^{-\alpha_1}}.
\end{eqnarray*}

Irreducible highest weight modules:
Let $\lambda=-n\lambda_1$ for some $n \in \bbZ^+$.  Since $\lambda$ is in the 
Wallach region, taking $n=0$ in the above formula:
$$ch_s L( \lambda ) = ch_s M( \lambda ) = \frac{e^{\lambda-\rho}}
{1+e^{-\alpha}}.$$
According to Theorem \ref{signedKL_formula},
$$1 = P^{\lambda,w_0}_{w_0,w_0} = \sgn (-w_0 \rho, \alpha_1 ) \varepsilon(
H_{\alpha_1,n}, s_1 ) P^{\lambda,w}_{w_0s_1,w_0} = \delta_{\alpha_1}^n
P^{\lambda,w}_{w_0s_1,w_0} = P^{\lambda,w_0}_{w_0s_1,w_0}$$
by Lemma 5.2.17 or Theorem 6.12 of \cite{Y}.  Substituting the values we 
have computed into Theorem \ref{hwm_formula}:
\begin{eqnarray*}
ch_s L( s_1 \lambda ) &=& R^{A(s_1\lambda,w_0)}( s_1\lambda) - P^{\lambda,w}
_{w_0s_1,w_0}R^{A(\lambda,w_0)}(\lambda) \\
&=& R^{A(n\lambda_1, w_0)}(s_1\lambda) - R^{A(-n\lambda_1,w_0)}(s_1 \lambda
- n \alpha_1 ) \\
&=& R^{A(n\lambda_1,w_0)}(s_1 \lambda ) - R^{A(0\lambda_1, w_0 )}(s_1 \lambda 
-n\alpha_1 ) \\
&=& \left( e^{s_1\lambda-\rho} + \cdots + e^{s_1 \lambda - (n-1)\alpha_1 - 
\rho } + \frac{e^{s_1\lambda-n\alpha_1-\rho}}{1+e^{-\alpha_1}} \right) -
\left( \frac{e^{s_1\lambda - n\alpha_1 -\rho}}{1+e^{-\alpha_1}} \right) \\
&=& e^{s_1\lambda -\rho} + e^{s_1\lambda-\alpha_1-\rho} + \cdots + 
e^{s_1\lambda-(n-1)\alpha_1 -\rho} .
\end{eqnarray*}

\begin{figure}[h]
\setlength{\unitlength}{1mm}%

\begingroup\makeatletter\ifx\SetFigFont\undefined%
\gdef\SetFigFont#1#2#3#4#5{%
  \reset@font\fontsize{#1}{#2pt}%
  \fontfamily{#3}\fontseries{#4}\fontshape{#5}%
  \selectfont}%
\fi\endgroup%
\begin{picture}(120,70)(5,0)
{\color[rgb]{0,0,0}\put(25,43){\line( 1,0){80}}
}%
{\color[rgb]{0,0,0}\multiput(15,43)(1,0){10}{\line( 1,0){.5}}
}%
{\color[rgb]{0,0,0}\multiput(105,43)(1,0){10}{\line( 1,0){.5}}
}%
\put(10,65){\makebox(0,0)[lb]{\smash{{\SetFigFont{10}{14.4}{\rmdefault}
{\mddefault}{\updefault}{\color[rgb]{0,0,0}$
\begin{array}[t]{cccccccc}
 & e^\mu & & e^\mu & & e^\mu & & \\
 &       & & +e^{\mu - \alpha_1} & & +e^{\mu - \alpha_1} & & \\
 &       & & & & +e^{\mu - 2 \alpha_1} & & \cdots \\
 & \lambda_1 & & 2\lambda_1 & & 3\lambda_1 & & \\ 
 & \bullet & & \bullet & & \bullet & & \\ 
 & H_{\alpha_1,1} & & H_{\alpha_1,2} & & H_{\alpha_1,3} & & \\
e^\mu & & e^\mu & & e^\mu & & e^\mu & \cdots \\
- e^{\mu-\alpha_1} & & + e^{\mu-\alpha_1} & & + e^{\mu-\alpha_1} & &
+ e^{\mu-\alpha_1} &  \\
+ e^{\mu-2\alpha_1} & & - e^{\mu-2\alpha_1} & & + e^{\mu-2\alpha_1} & &
+ e^{\mu-\alpha_1} &  \\
- e^{\mu-3\alpha_1} & & + e^{\mu-3\alpha_1} & & - e^{\mu-3\alpha_1} & &
+ e^{\mu-3\alpha_1} &  \\
+ e^{\mu-4\alpha_1} & & - e^{\mu-4\alpha_1} & & + e^{\mu-4\alpha_1} & &
- e^{\mu-4\alpha_1} &  \\
\vdots & & \vdots & & \vdots & & \vdots & \\
\end{array}
$}%
}}}}
\put(5,58){\makebox(0,0)[lb]{\smash{{\SetFigFont{10}{14.4}{\rmdefault}
{\mddefault}{\updefault}{\color[rgb]{0,0,0}$ch_s L(\mu)$}%
}}}}
\put(5,38){\makebox(0,0)[lb]{\smash{{\SetFigFont{10}{14.4}{\rmdefault}
{\mddefault}{\updefault}{\color[rgb]{0,0,0}$ch_s M(\mu)$}%
}}}}
\put(117,43){\makebox(0,0)[lb]{\smash{{\SetFigFont{10}{14.4}{\rmdefault}
{\mddefault}{\updefault}{\color[rgb]{0,0,0}$\mu$}%
}}}}
\end{picture}%
\caption{$\mathfrak{su}(2)$}
\end{figure}
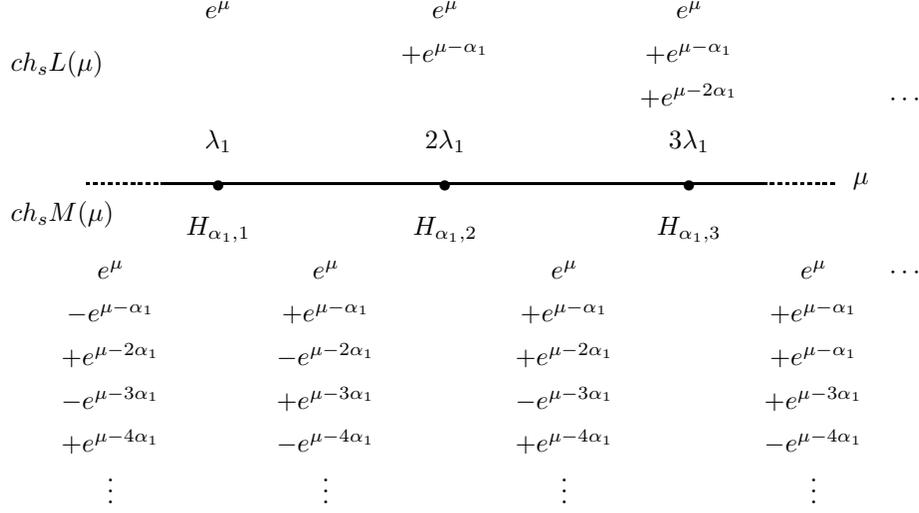

\subsection*{Example 2:  $\g_0 = \mathfrak{sl}(2,\bbR)$.}
We may proceed as in the previous example, but substitute $\delta_{\alpha_1}
= -1$ instead of $\delta_{\alpha_1} = 1$.

Irreducible Verma modules:  For $\lambda \in \h^*$ such that $( \lambda, 
\alpha_1^\vee) \in (n,n+1)$ where $n \in \bbZ_{\geq 0}$:
\begin{eqnarray*}
ch_s M( \lambda ) &=& R^{A(n\lambda_1, w_0)}( \lambda ) = 
\frac{\sum_{i=1}^n (-1)^{n-i+1} 2 e^{\overline{r_i}r_i\lambda-\rho} + 
e^{\lambda-\rho}}
{\displaystyle{\prod_{\alpha \in \Delta^+( \mfp, \mft)}( 1 - e^{-\alpha})
\prod_{\alpha \in \Delta^+( \mfk, \mft)}( 1 + e^{-\alpha})}} \\
&=& 
\frac{\sum_{i=1}^n (-1)^i 2 e^{\lambda - i \alpha_1 -\rho} + e^{\lambda-\rho}}
{1-e^{-\alpha_1}} \\
&=& e^{\lambda-\rho} - e^{\lambda-\rho-\alpha_1 -\rho} + \cdots +
(-1)^{n-1}e^{\lambda-(n-1)\alpha_1 -\rho} + (-1)^n \frac{e^{\lambda-n\alpha_1
-\rho}} {1-e^{-\alpha_1}}.
\end{eqnarray*}

Irreducible highest weight modules:  For $\lambda=-n\lambda_1$ where $n \in
\bbZ^+$:
$$ ch_s L( \lambda ) = ch_s M( \lambda ) = \frac{e^{\lambda-\rho}}
{1-e^{-\alpha_1}}.$$
Since $P^{\lambda,w_0}_{w_0s_1,s_0} = (-1)^n$, we have
\begin{eqnarray*}
ch_s L( s_1 \lambda) &=& R^{A(s_1\lambda,w_0)}(s_1\lambda) - P^{\lambda,w_0}
_{w_0s_1, w_0}R^{A(\lambda,w_0)}(\lambda) \\
&=& \left( \sum_{i=0}^{n-1} (-1)^i e^{s_1 \lambda - i\alpha_1 -\rho} + 
(-1)^n \frac{e^{s_1\lambda-n\alpha_1 -\rho}}{1-e^{-\alpha_1}} \right)
-(-1)^n \left( \frac{e^{s_1\lambda-n \alpha_1 - \rho}}{1 - e^{-\alpha_1}}
\right) \\
&=& e^{s_1\lambda - \rho} - e^{s_1\lambda-\alpha_1 -\rho} + \cdots + 
(-1)^{n-1}e^{s_1\lambda - (n-1)\alpha_1 - \rho}.
\end{eqnarray*}
\begin{figure}[h]
\setlength{\unitlength}{1mm}%

\begingroup\makeatletter\ifx\SetFigFont\undefined%
\gdef\SetFigFont#1#2#3#4#5{%
  \reset@font\fontsize{#1}{#2pt}%
  \fontfamily{#3}\fontseries{#4}\fontshape{#5}%
  \selectfont}%
\fi\endgroup%
\begin{picture}(120,70)(5,0)
{\color[rgb]{0,0,0}\put(25,43){\line( 1,0){80}}
}%
{\color[rgb]{0,0,0}\multiput(15,43)(1,0){10}{\line( 1,0){.5}}
}%
{\color[rgb]{0,0,0}\multiput(105,43)(1,0){10}{\line( 1,0){.5}}
}%
\put(10,65){\makebox(0,0)[lb]{\smash{{\SetFigFont{10}{14.4}{\rmdefault}
{\mddefault}{\updefault}{\color[rgb]{0,0,0}$
\begin{array}[t]{cccccccc}
 & e^\mu & & e^\mu & & e^\mu & & \\
 &       & & -e^{\mu - \alpha_1} & & -e^{\mu - \alpha_1} & & \\
 &       & & & & +e^{\mu - 2 \alpha_1} & & \cdots \\
 & \lambda_1 & & 2\lambda_1 & & 3\lambda_1 & & \\ 
 & \bullet & & \bullet & & \bullet & & \\ 
 & H_{\alpha_1,1} & & H_{\alpha_1,2} & & H_{\alpha_1,3} & & \\
e^\mu & & e^\mu & & e^\mu & & e^\mu & \cdots \\
+ e^{\mu-\alpha_1} & & + e^{\mu-\alpha_1} & & + e^{\mu-\alpha_1} & &
+ e^{\mu-\alpha_1} &  \\
+ e^{\mu-2\alpha_1} & & - e^{\mu-2\alpha_1} & & + e^{\mu-2\alpha_1} & &
+ e^{\mu-\alpha_1} &  \\
+ e^{\mu-3\alpha_1} & & - e^{\mu-3\alpha_1} & & - e^{\mu-3\alpha_1} & &
+ e^{\mu-3\alpha_1} &  \\
+ e^{\mu-4\alpha_1} & & - e^{\mu-4\alpha_1} & & - e^{\mu-4\alpha_1} & &
- e^{\mu-4\alpha_1} &  \\
\vdots & & \vdots & & \vdots & & \vdots & \\
\end{array}
$}%
}}}}
\put(5,58){\makebox(0,0)[lb]{\smash{{\SetFigFont{10}{14.4}{\rmdefault}
{\mddefault}{\updefault}{\color[rgb]{0,0,0}$ch_s L(\mu)$}%
}}}}
\put(5,38){\makebox(0,0)[lb]{\smash{{\SetFigFont{10}{14.4}{\rmdefault}
{\mddefault}{\updefault}{\color[rgb]{0,0,0}$ch_s M(\mu)$}%
}}}}
\put(117,43){\makebox(0,0)[lb]{\smash{{\SetFigFont{10}{14.4}{\rmdefault}
{\mddefault}{\updefault}{\color[rgb]{0,0,0}$\mu$}%
}}}}
\end{picture}%
\caption{$\mathfrak{sl}(2)$}
\end{figure}
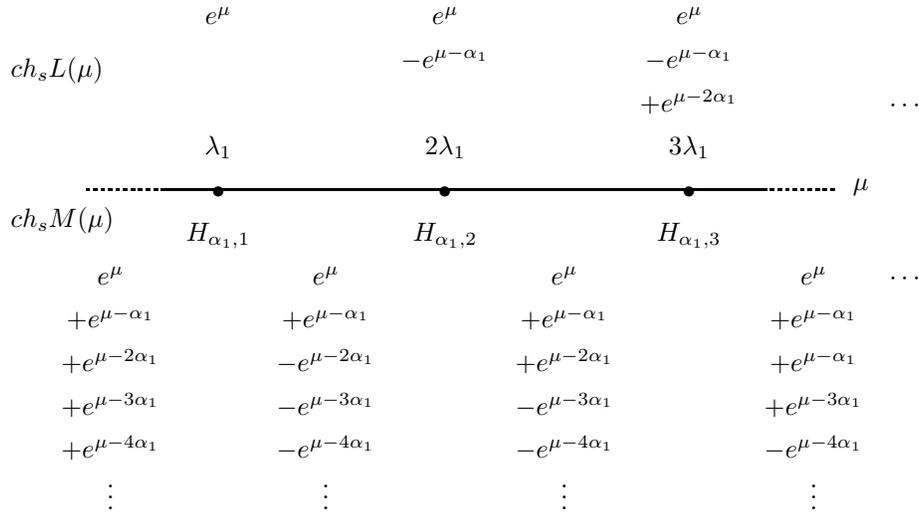
\nocite{V}
\nocite{KV}
\nocite{H}
\nocite{H2}
\nocite{H3}
\nocite{V2}
\nocite{V4}
\nocite{D}
\nocite{K}
\nocite{J3}
\nocite{SV}
\nocite{MFF}
\nocite{EW}
\bibliographystyle{alpha}
\bibliography{hwm_signature}

\begin{thebibliography}{Vog79b}

\bibitem[Bar83]{BA}
Dan Barbasch.
\newblock Filtrations on {V}erma modules.
\newblock {\em Ann. Sci. \'Ecole Norm. Sup. (4)}, 16(3):489--494, 1983.

\bibitem[BB81]{BB2}
Alexandre Be{\u\i}linson and Joseph Bernstein.
\newblock Localisation de {$g$}-modules.
\newblock {\em C. R. Acad. Sci. Paris S\'er. I Math.}, 292(1):15--18, 1981.

\bibitem[BB93]{BB}
A.~Beilinson and J.~Bernstein.
\newblock {\em A proof of {J}antzen conjectures}, volume~16 of {\em Adv. Soviet
  Math.}
\newblock Amer. Math. Soc., Providence, RI, 1993.

\bibitem[BGG71]{BGG2}
I.~N. Bernstein, I.~M. Gelfand, and S.~I. Gelfand.
\newblock Structure of representations that are generated by vectors of higher
  weight.
\newblock {\em Funckcional. Anal. i Prilo\v zen.}, 5(1):1--9, 1971.

\bibitem[BGG76]{BGG}
I.~N. Bernstein, I.~M. Gelfand, and S.~I. Gelfand.
\newblock A certain category of $\g$-modules.
\newblock {\em Functional Anal. Appl.}, 10:87--92, 1976.

\bibitem[BK81]{KB}
J.-L. Brylinski and M.~Kashiwara.
\newblock Kazhdan-{L}usztig conjecture and holonomic systems.
\newblock {\em Invent. Math.}, 64(3):387--410, 1981.

\bibitem[Dix96]{D}
Jacques Dixmier.
\newblock {\em Enveloping Algebras}.
\newblock Number~11 in Graduate Studies in Mathematics. American Mathematical
  Society, Providence, RI, 1996.

\bibitem[DL77]{DL}
Vinay Deodhar and James Lepowsky.
\newblock On multiplicity in the {J}ordan-{H}\"{o}lder series of {V}erma
  modules.
\newblock {\em J. Algebra}, 49(2):512--524, 1977.

\bibitem[EW80]{EW}
Thomas Enright and Nolan~R. Wallach.
\newblock Notes on homological algebra and representations of {L}ie algebras.
\newblock {\em Duke Math. J.}, 47(1):1--15, 1980.

\bibitem[GJ81]{GJ}
O.~Gabber and A.~Joseph.
\newblock Towards the {K}azhdan-{L}usztig conjecture.
\newblock {\em Ann. Sci \'{E}cole Norm. Sup. (4)}, 14(3):261--302, 1981.

\bibitem[Hum72]{H3}
James~E.\ Humphreys.
\newblock {\em Introduction to Lie Algebras and Representation Theory}.
\newblock Number~9 in Graduate Texts in Mathematics. Springer-Verlag, New York,
  1972.

\bibitem[Hum78]{H2}
James~E. Humphreys.
\newblock Finite and infinite dimensional modules for semisimple {L}ie
  algebras.
\newblock In {\em Lie theories and their applications (Proc. Ann. Sem. Canad.
  Math. Congr., Queen's Univ., Kingston, Ont., 1977)}, volume~48 of {\em
  Queen's Papers in Pure and Appl. Math.}, pages 1--64. Queen's Univ.,
  Kingston, Ont., 1978.

\bibitem[Hum90]{H}
James~E.\ Humphreys.
\newblock {\em Reflection groups and {C}oxeter groups}.
\newblock Number~29 in Cambridge Studies in Advanced Mathematics. Cambridge
  University Press, Cambridge, 1990.

\bibitem[Jan73]{J3}
Jens~Carsten Jantzen.
\newblock Darstellungen halbeinfacher algebraischer {G}ruppen und zugeordnete
  kontravariante {F}ormen.
\newblock {\em Bonn. Math. Schr.}, (67):v+124, 1973.

\bibitem[Jan74]{J2}
Jens~Carsten Jantzen.
\newblock Zur {C}harakterformel gewisser {D}arstellungen halbeinfacher
  {G}ruppen und {L}ie-{A}lgebren.
\newblock {\em Math. Z.}, 140:127--149, 1974.

\bibitem[Jan79]{J}
Jens~Carsten Jantzen.
\newblock {\em Moduln mit einem h\"{o}chsten {G}ewicht}.
\newblock Number 750 in Lecture Notes in Mathematics. Springer, Berlin, 1979.

\bibitem[KL79]{KL}
D.~Kazhdan and G.~Lusztig.
\newblock Representations of {C}oxeter groups and {H}ecke algebras.
\newblock {\em Invent. Math.}, 53:165--184, 1979.

\bibitem[Kna96]{K}
Anthony~W. Knapp.
\newblock {\em Lie Groups Beyond an Introduction}.
\newblock Number 140 in Progress in Mathematics. Birkh\"auser, Boston, 1996.

\bibitem[Kos75]{KO}
Bertram Kostant.
\newblock On the tensor product of a finite and an infinite dimensional
  representation.
\newblock {\em J. Functional Analysis}, 20(4):257--285, 1975.

\bibitem[KV95]{KV}
Anthony~W. Knapp and David~A. Vogan.
\newblock {\em Cohomological induction and unitary representations}.
\newblock Number~45 in Princeton Mathematical Series. Princeton University
  Press, Princeton, New Jersey, 1995.

\bibitem[MFF86]{MFF}
F.G. Malikov, B.L. Feigin, and D.B. Fuks.
\newblock Singular vectors in {Verma} modules over {Kac-Moody} algebras.
\newblock {\em Functional Analysis Applications}, 20(2):103--113, 1986.

\bibitem[SV80]{SV}
Birgit Speh and David~A. Vogan.
\newblock Reducibility of generalized principal series representations.
\newblock {\em Acta Math.}, 145(3-4):227--299, 1980.

\bibitem[Ver68]{VE}
Daya-Nand Verma.
\newblock Structure of certain induced representations of complex semisimple
  {L}ie algebras.
\newblock {\em Bull. Amer. Math. Soc.}, 74:160--166, 1968.

\bibitem[Vog79a]{V2}
David~A. Vogan.
\newblock Irreducible characters of semisimple {L}ie groups {I}.
\newblock {\em Duke Math. J.}, 46:61--108, 1979.

\bibitem[Vog79b]{V4}
David~A. Vogan.
\newblock Irreducible characters of semisimple {L}ie groups {II}. {T}he
  {K}azhdan-{L}usztig conjectures.
\newblock {\em Duke Math. J.}, 46:805--859, 1979.

\bibitem[Vog81]{V}
David~A. Vogan.
\newblock {\em Representations of real reductive {L}ie groups}.
\newblock Number~15 in Progress in Mathematics. Birkh\"auser, Boston,
  Massachusetts, 1981.

\bibitem[Vog84]{V3}
David~A. Vogan.
\newblock Unitarizability of certain series of representations.
\newblock {\em Ann. of Math. (2)}, 120(1):141--187, 1984.

\bibitem[Wal84]{W}
Nolan~R. Wallach.
\newblock On the unitarizability of derived functor modules.
\newblock {\em Invent. Math.}, 78(1):131--141, 1984.

\bibitem[Yee05]{Y}
Wai~Ling Yee.
\newblock The signature of the {S}hapovalov form on irreducible {V}erma
  modules.
\newblock {\em Representation Theory}, 9:638--677, 2005.

\end{thebibliography}
\end{document}